\documentclass[pdftex]{amsart}
\usepackage{amsmath}
\usepackage{amssymb}
\usepackage{booktabs}
\usepackage{ifthen}
\newboolean{pdfOutput}
\setboolean{pdfOutput}{true}
\newboolean{abridged}
\setboolean{abridged}{true}
\usepackage{euscript}
\usepackage[all]{xy}
\usepackage{varioref}
\usepackage{rotating}
\newtheorem{theorem}{Theorem}
\newtheorem{corollary}{Corollary}
\newtheorem{lemma}{Lemma}
\newtheorem{proposition}{Proposition}

\newtheorem{definition}{Definition}
\newtheorem{remark}{Remark}
\newcommand{\openSet}{\textrm{ open }}
\newcommand{\pd}[2]{\ensuremath{ \frac{\partial {#1}}{\partial {#2}} }}
\newcommand{\R}[1]{\ensuremath{\mathbb{R}^{#1}}}
\newcommand{\C}[1]{\ensuremath{\mathbb{C}^{#1}}}
\newcommand{\CP}[1]{\ensuremath{\mathbb{CP}^{#1}}}
\newcommand{\Gro}[2]{\widetilde{\operatorname{Gr}}\left ( #1, #2 \right ) }
\newcommand{\GrC}[2]{\operatorname{Gr}_{\C{}} \left ( #1, #2 \right )}
\newcommand{\GL}[1]{ \operatorname{GL} \left ( {#1} \right ) }

\newcommand{\Lin}[2]{ \operatorname{Lin} \left ( {#1} , {#2} \right ) }

\newcommand{\nForms}[2]{\ensuremath{\Omega^{#1} \left ( {#2} \right )}}
\newcommand{\Cohom}[2]{\ensuremath{
H^{#1} \left ( {#2} \right )}}

\newcommand{\trans}[1]{{}^t{#1}}
\newcommand{\bari}{\bar{\imath}}
\newcommand{\barj}{\bar{\jmath}}
\newcommand{\bark}{\bar{k}}
\newcommand{\barm}{\bar{m}}
\newcommand{\barp}{\bar{p}}
\newcommand{\barmu}{\bar{\mu}}
\newcommand{\barnu}{\bar{\nu}}
\newcommand{\barsigma}{\bar{\sigma}}
\newcommand{\bartau}{\bar{\tau}}
\newcommand{\barepsilon}{\bar{\epsilon}}
\DeclareMathOperator{\Det}{Det}

\DeclareMathOperator{\Vertical}{Vert}
\DeclareMathOperator{\Real}{Re}
\DeclareMathOperator{\Imag}{Im}

\newcommand{\extraSection}[2]
{
\ifthenelse{\boolean{abridged}}
  {
  }
  {
    \section{#1} 
    \begin{center}
    \emph{This section will not be referred to
    subsequently, and may be skipped.}
    \end{center}
    \par{#2}
  }
}

\newcommand{\extraSubsection}[2]
{
\ifthenelse{\boolean{abridged}}
  {
  }
  {
    \subsection{#1}
    \begin{center}
    \emph{This subsection will not be referred to
    subsequently, and may be skipped.}
    \end{center}
    \par{#2}
  }
}

\newcommand{\extraStuff}[1]%
{%
\ifthenelse{\boolean{abridged}}%
  {%
  }%
  {%
    {#1}%
  }%
}

\begin{document}
\title{Analogues of Complex Geometry}
\author{Benjamin McKay}
\address{University College Cork}
\email{b.mckay@ucc.ie}
\date{\today}
\begin{abstract}
We prove that there are no 
pseudoholomorphic theories of anything
other than curves, even if one allows more
general spaces than almost complex manifolds.
The proof is elementary, except for theories
of pseudoholomorphic hypersurfaces, where
topological techniques are needed.
Surprisingly, hypersurface theories
exist ``microlocally'' (in great abundance) 
to all orders perturbatively, but not ``locally.''
\end{abstract}
\maketitle
\tableofcontents
\section{Introduction}
Misha Gromov, asked to point to the future directions
of geometry, responded \cite{Gromov:2000} that 
the geometry of solutions of first-order
systems of partial differential equations, generalizing
(1) the Cauchy--Riemann equations of almost complex manifolds  
and (2) the equations of calibrated submanifolds in
exceptional holonomy manifolds, is one of the
key directions to follow. This paper is the first
attempt since Gromov's paper to follow that direction.
\subsection{The problem}
This paper is a study of the system of
differential equations
\begin{equation} \label{eqn:Eqn}
\pd{w^i}{z^{\bar{\mu}}}
=
F^i_{\bar{\mu}}
\left ( 
z^{\nu}, z^{\bar{\nu}},
w^j, w^{\barj},
\pd{w^j}{z^{\nu}},
\pd{w^{\barj}}{z^{\bar{\nu}}}
\right ) \tag{$\blacklozenge$}
\end{equation}
where $(z,w) \in \openSet \subset \C{n+d}$ and
$w \left ( z,\bar{z} \right )$  is a complex
valued function of the complex $z$ variables.
\extraStuff{%
As a mnemonic, the number of 
dependent variables is $d$, and the number of
independent variables is $n$.}
We will assume that the function 
\[
F : \openSet \subset \C{n} \times\C{d} \times \C{nd} \to \C{nd}
\]
is smooth enough to carry through our 
arguments; $F$ four times continuously differentiable
will suffice.

If $F=0$, then these are the
Cauchy--Riemann equations of complex analysis.
Thus we are studying the deformation theory 
of the Cauchy--Riemann
equations through elliptic equations.
The idea is to try to liberate the 
Cauchy--Riemann equations from the
global rigidity problems associated with
complex structures, producing a more
flexible theory, almost topological, holding onto only
the most fundamental analytical estimates
of the Cauchy--Riemann equations.

The problem solved in this paper is to find
the systems of such equations which have 
the same \emph{tableau}
as the Cauchy--Riemann equations (in the sense of 
Cartan, see Cartan \cite{Cartan:1945}). This is 
a purely formal requirement, i.e. it requires
only the vanishing of certain 
algebraic expressions in $F$ and its 
first derivatives. Another way to look at it:
we will show that equations~\vref{eqn:Eqn} 
have Cauchy--Riemann tableau at a point precisely
it is possible to change the coordinates
to have the function $F$ vanish along with
its first derivatives at that point. 

The requirement of Cauchy--Riemann tableau
is natural: if $n=1$ or $d=1$ then an
\emph{a priori} bound on the mean curvature 
of 2-jets of solutions of \ref{eqn:Eqn} 
in any Riemannian metric is equivalent to having 
Cauchy--Riemann tableau. This follows
from the observation that the mean curvature
bound gets smaller as we dilate, while the
tableau is invariant under dilation, and 
then following
J.M.~Landsberg's paper \cite{Landsberg:1992}.
In that paper, he demonstrated that the tableaux
for $2dn$ equations for $2d$ functions
of $2n$ variables whose integral manifolds 
are minimal in a flat metric are the
Cauchy--Riemann tableaux, when $n=1$ or $d=1$.\footnote{In 
all other cases,
$n>1$ and $d>1$, it is not known if mean
curvature bounds are equivalent to Cauchy--Riemann
tableau, but would follow from conjectures of
Landsberg.} Mean curvature bounds are
vital to all known proofs of all 
Gromov-type compactness theorems.

Overdetermined systems of differential equations
rarely have any solutions, even locally. One
generally obtains obstructions to the existence
of solutions by examining the \emph{torsion}
of the tableau. However,
the torsion of the Cauchy--Riemann equations vanishes.
Therefore so does the torsion of any equations with the same
tableau. It follows (see Cartan \cite{Cartan:1945})
that systems of equations with Cauchy--Riemann tableau
have no finite order
obstructions to local solvability.
Unobstructed solvability 
makes their study even more compelling.
\subsection{The solution}
\begin{theorem}
Every system of partial 
differential equations of the form 
\ref{eqn:Eqn} for at least two complex functions
$w$ of at least two complex variables $z$ 
which has Cauchy--Riemann tableau
becomes the Cauchy--Riemann equations
in some system of local coordinates.
\end{theorem}
Strangely, when there is only one $w$ or
one $z$ variable, then there are
many equations with Cauchy--Riemann
tableau. Call
a set of first-order equations 
\emph{proper} if 
the 1-jets satisfying
them at any chosen point $(z,w)$ (i.e.
fixed 0-jet) form a compact set.
The Cauchy--Riemann equations are proper.
\begin{theorem}
Proper equations
with Cauchy--Riemann tableau
for functions
of several variables
become Cauchy--Riemann equations
after a coordinate change. 
\end{theorem}
The proof uses the theory of complex contact
geometry (due to Merkulov), computes
Chern numbers of some complex vector bundles,
and uses a theorem of the author classifying
the diffeomorphism types of great
circle fibrations of spheres.
Consequently, the only 
possibility left, where one can
globally perturb the Cauchy--Riemann
equations, fixing their tableau, is that of equations 
with one independent variable, i.e.
theories of pseudoholomorphic curves. 
Following Landsberg, pseudoholomorphic
curve theories are our only hope for
Gromov-type compactness theorems.
\subsection{The big picture}
These theorems limit the possibilities
for constructing analogues of complex geometry,
and (for the sake of symplectic topology) constructing
analogues of K\"ahler geometry. Researchers
will naturally feel encouraged by these theorems
to study pseudoholomorphic curves, following
Gromov.

Note the curious
appearance of equations in dimension 4
(i.e. when there is only one $z$ and one $w$ variable,
in local coordinates),
generalizing pseudoholomorphic curves. These were probably 
first noticed by Lavrentiev \cite{Lavrentiev:1947,Lavrentiev:1948}.
Smooth topological projective planes provide 
one plentiful source of these equations; see McKay \cite{McKay:2005}.
Clearly these
equations in dimension 4 demand investigation.

Boris Kruglikov has similar results for almost complex manifolds
(see \cite{Kruglikov:2003} p. 66 for a comparison with this paper).
Our more remarkable results: (1)
the equations arising in dimension 4 and (2)
the microlocal presence but local absence
of hypersurface equations, do not arise
in Kruglikov's paper. These results do not
concern almost complex structures (which
have only quasi-linear Cauchy--Riemann equations), but
fully nonlinear systems of differential
equations.

\section{Constructing $G$-structures on the Grassmann bundles} 
Given manifolds $M$ and $N$, and a smooth map
$f : M \to N$, write $f'(m) : T_m M \to T_{f(m)} N$ 
for the induced map on tangent spaces.
Consider a manifold $M$ of dimension $n+d$, and
the bundle 
\[
\pi : \Gro{n}{TM} \to M
\]   
whose fiber above a point $m \in M$ 
is the Grassmann bundle $\Gro{n}{T_m M}$ of $n$ dimensional
oriented vector subspaces (which we will call $n$-planes) in the
tangent spaces of $M$. 
\extraStuff{
The dimensions are
\[
\xymatrix{
\Gro{n}{T_m M}^{nd} \ar[r] & \Gro{n}{TM}^{nd+n+d} \ar[d]^{\pi} \\
& M^{n+d}
}
\]
}
There is a canonical field of $n+nd$-planes on
$\Gro{n}{TM}$ given by
\[
\Theta_P = \pi'(P)^{-1} P
\] 
called the \emph{polycontact plane field}
and a canonical isomorphism
\[
\pi'(P) : \Theta_P / \ker \pi'(P) \to P.
\]
If we have any immersion $f : \Sigma \to M$
of a manifold $\Sigma$ of dimension $n$,
then there is a canonical lift $\hat{f} : \Sigma \to \Gro{n}{TM}$
defined by
\[
\hat{f}(s) = f'(s) \, T_s \Sigma
\]
and it is clear that 
\[
\hat{f}'(s) : T_s \Sigma \to \Theta_{\hat{f}(s)}.
\]
Conversely, if $F : \Sigma \to \Gro{n}{TM}$ is
any immersion satisfying 
\[
F'(s) : T_s \Sigma \to \Theta_{F(s)}
\]
and $F$ is transverse to
the fibers of $\Gro{n}{TM} \to M$, then
\[
F=\hat{f}
\]
where $f = \pi F$.

We also know that the tangent spaces of the fibers 
have canonical identifications
\[
T_P \Gro{n}{T_m M} \cong \Lin{P}{T_m M/P}
\]
with the spaces of linear maps from
$P$ to $T_m M/P$
given in the following manner.
Let $P(t)$ be any family of $n$-planes
in $T_m M$, and $\phi(t) : T_m M \to W$ any
family of linear maps so that
\[
\ker \phi(t) = P(t).
\]
Then let 
\[
[\phi](t) : v+ P(t) \in T_mM/P(t) \to \phi(t)(v) \in W. 
\]
Identify $P'(t)$ with
\[
[\phi](t)^{-1} \left . \phi'(t) \right |_P : P \to T_m M/P.
\]
\extraStuff{
We have a diagram
\[
\xymatrix{
& 0 \ar[d] & 0 \ar[d] & 0 \ar[d] & \\
0 \ar[r] & \Lin{P}{T_m M/P} \ar[r] \ar[d] & \Theta_P \ar[r] \ar[d] & 
P \ar[r] \ar[d] & 0 \\
0 \ar[r] & T_P \Gro{n}{T_m M} \ar[r] \ar[d] & T_P \Gro{n}{TM} \ar[r] \ar[d] & 
T_m M  \ar[r] \ar[d] & 0 \\
  & 0 \ar[r] & T_P \Gro{n}{TM} / \Theta_P \ar[r] \ar[d] &
T_m M / P \ar[r] \ar[d] & 0 \\
 & & 0 & 0 & \\
}
\]
}

Let $V=\R{n+d}$ and $P_0 =\R{n} \subset V$.
Any choice of linear isomorphism
\[
u : T_m M \to V
\]
taking an $n$-plane $P \subset T_m M$ to $u(P) = P_0$
induces an isomorphism on the tangent space of the fiber
\[
u' : T_P \Gro{n}{T_m M} = \Lin{P}{T_m M/P} \to \Lin{P_0}{V/P_0}.
\] 
given by
\[
[\phi](t)^{-1} \phi'(t) \mapsto
\left [ \phi(t) u^{-1} \right ]^{-1} \phi'(t) u^{-1}.
\]
Let $H \subset \GL{V}$ be the 
subgroup of linear transformations
leaving the plane $P_0$ invariant.
If we change the choice of isomorphism
$u$ to another, say $v$, which still
identifies $P$ with $P_0 \subset V$, then
\[
v = g u
\]
where $g \in H$.
We have the obvious homomorphisms
\[
\rho_{P_0} : H \to \GL{V}
\]
and 
\[
\rho_{V/P_0} : H \to \GL{V/P_0}
\]
and we find that
\[
\left ( gu \right )' = \rho_{V/P_0}(g) u' \rho_{P_0}^{-1}(g).
\]

Given such a map
\[
u : T_m M \to V
\] 
we will consider an \emph{adapted 
coframe} on the total space $\Gro{n}{TM}$
to be a linear isomorphism
\[
U : T_P \Gro{n}{TM} \to V \oplus \Lin{P_0}{V/P_0}
\]
so that the $V$ part vanishes
on the fibers, and hence is defined on the
base $T_m M$, and so that moreover it
equals $u$ on the base, and so that the
$\Lin{P_0}{V/P_0}$ part
equals $u'$ on the fiber. 
The part valued in $\Lin{P_0}{V/P_0}$ 
is not determined completely by this
condition.
Such a map $U$ is determined by $u$ up to
choices of adding some linear function
\[
T_P \Gro{n}{TM} \to \Lin{P_0}{V/P_0}
\] 
on the tangent space which vanishes on
the fiber, and consequently defined
on the base $T_mM$. The base is identified
by $u$ with $V$, so 
an adapted coframe $U$ is uniquely
determined by the map $u$ up to choice of 
a linear map
\[
\Lin{V}{\Lin{P_0}{V/P_0}}.
\]
Now if we change $u$, so that we
pick some other map 
\[
v : T_m M \to V
\]
which identifies the same $n$-plane $P \subset T_m M$
with $P_0 \subset V$,
then we must have 
\[
v = gu
\]
for some $g \in H$.
What is the effect on the map $U$?
Its $V$ valued part is changed
by $g$, and its $\Lin{P_0}{V/P_0}$
is changed by action of
\[
\rho_{V/P_0}(g) \otimes \trans{\rho_{P_0}}^{-1}.
\]
Moreover we will still have an adapted
coframe if we alter this one by plugging
the $V$ valued part into any element
of $\Lin{V}{\Lin{P_0}{V/P_0}}$ and
adding this to the $\Lin{P_0}{V/P_0}$
part. Hence the adapted coframes are
well defined up to this action of the
group \( G_0 = H \ltimes \Lin{V \otimes P_0}{V/P_0}. \)
Let $B_0$ be the bundle of all adapted 
coframes and 
\[
\Pi : B_0 \to \Gro{n}{TM}
\]
be the obvious map.
We have explained how to build a left action
of $G_0$ on $B_0$.
Henceforth we will instead 
let $G_0$ act on $B$ on the right
by using the inverse of the left action:
\[
r_g U = g^{-1} U.
\]
We have therefore found that
$B_0$ is a principal right $G_0$ bundle.
We define the \emph{soldering 1-form}
\[
\omega \in \nForms{1}{B} \otimes
\left ( V \oplus \Lin{P_0}{V/P_0} \right )
\]
by the equation
\[
\omega_U = U \Pi'.
\]
\extraStuff{
It satisfies 
\[
r_g^* \omega = g^{-1} \omega.
\]
}

Roughly put, because the soldering form is invariantly defined,
any equations that we can write in terms of the soldering
form are necessarily invariant under diffeomorphism.
The soldering form provides a diffeomorphism invariant
computational apparatus. To relate it to our differential
equations \ref{eqn:Eqn} we need to complete the tedious
task of expressing the soldering form
and its exterior derivative in local coordinates.

In local coordinates 
\[
x^{\mu}, y^i
\]
on $M$ near a point $m \in M$, 
we find that every $n$ plane on
which the $dx^{\mu}$ are independent 1-forms
has the form
\[
dy^i = p^i_{\mu} dx^{\mu}.
\]
Moreover these $p^i_{\mu}$ are arbitrary,
and therefore provide local coordinates
on $\Gro{n}{T_m M}$. Consequently
\[
x^{\mu}, y^i, p^i_{\mu}
\]
are coordinates on $\Gro{n}{TM}$
near $P$ where 
\[
P = \left \{ dy^i = 0 \right \} \subset T_m M.
\]
Take a basis $e_{\mu}$ of $P_0$ and extend to
a basis $e_{\mu},e_i$ of $V$.
The 1-forms
\begin{align*}
\eta^{i} =& dy^i - p^i_{\mu} \, dx^{\mu} \\
\eta^{\mu} =& dx^{\mu} \\
\eta^i_{\mu} =& dp^i_{\mu} 
\end{align*}
put together\footnote{
The reader must remain on guard for confusion
in this notation. The Roman indices run
over $n+1,\dots,n+d$, while the Greek
run over $1,\dots,n$.}
form a local section $\eta$
of the bundle $B_0 \to \Gro{n}{TM}$,
with
\[
U(x,y,p) = \left ( \eta^{\mu} \oplus \eta^i
\oplus \eta^i_{\mu} \right )
: T_P \Gro{n}{TM} \to V \oplus \Lin{P_0}{V/P_0}.
\]
The bundle $B_0$ then admits local
coordinates
\[
x^{\mu}, y^i, p^i_{\mu}, a^i_j, a^{\mu}_j, a^{\mu}_{\nu}, 
a^i_{\mu j}, a^i_{\mu \nu}
\]
where we write every element of $B_0$ close
to our section as
\[
U = 
\begin{pmatrix}
a^i_j & 0 & 0 \\
a^{\mu}_j & a^{\mu}_{\nu} & 0 \\
a^i_{\mu j} & a^i_{\mu \nu} & a^{i \nu}_{\mu j}  
\end{pmatrix}
\begin{pmatrix}
\eta^j \\
\eta^{\nu} \\
\eta^j_{\nu}
\end{pmatrix}.
\]
The requirement that coframes transform
under the group $G_0$ forces
\[
 a^{i \nu}_{j \mu} =
 a^i_j A^{\nu}_{\mu}
\]
% \begin{align*}
% a^{i \nu}_{j \mu} =&
% a^i_j A^{\nu}_{\mu} \\
% a^i_{\mu \nu} a^{\mu}_{\sigma}
% =&
% a^i_{\mu \sigma} a^{\mu}_{\nu}
% \end{align*}
where we write $A^{\mu}_{\nu}$ for the inverse 
matrix of $a^{\mu}_{\nu}$:
\[
A^{\mu}_{\sigma} a^{\sigma}_{\nu} = \delta^{\mu}_{\nu}.
\]
Otherwise the $a$ are arbitrary, except that
the matrices $a^i_j$ and $a^{\mu}_{\nu}$ must
be invertible.
The soldering 1-form is given in
coordinates by the same expression:
\[
\begin{pmatrix}
\omega^i \\
\omega^{\mu} \\
\omega^i_{\mu} 
\end{pmatrix}
= 
\begin{pmatrix}
a^i_j & 0 & 0 \\
a^{\mu}_j & a^{\mu}_{\nu} & 0 \\
a^i_{\mu j} & a^i_{\mu \nu} & a^i_j A^{\nu}_{\mu}
\end{pmatrix}
\begin{pmatrix}
\eta^j \\
\eta^{\nu} \\
\eta^j_{\nu}
\end{pmatrix}
\]
or 
\[
\omega = a \eta
\]
for short. 
Write
\[
a^{-1} =
\begin{pmatrix}
A^j_k & 0 & 0 \\
A^{\nu}_k & A^{\nu}_{\sigma} & 0 \\
A^j_{\nu k} & A^j_{\nu \sigma} & A^j_k a^{\sigma}_{\nu} \end{pmatrix}.
\]
\extraStuff{
Note:
\begin{align*}
a^i_j A^j_k =& \delta^i_k \\
a^{\mu}_{\nu} A^{\nu}_{\sigma} =& \delta^{\mu}_{\sigma} \\
A^{\nu}_k =& - A^{\nu}_{\mu} a^{\mu}_j A^j_k \\
A^j_{\nu \sigma} =& 
- a^{\mu}_{\nu} A^j_i a^i_{\mu \tau} A^{\tau}_{\sigma} \\
A^j_{\nu k} =&
a^{\mu}_{\nu} A^j_m \left ( 
a^m_{\mu \sigma} A^{\sigma}_{\tau} a^{\tau}_i 
-
a^m_{\mu i}
\right ) A^i_k
\end{align*}
}
Differentiating, we find
\extraStuff{
\[
d \omega = da \wedge \eta + a \, d \eta 
\]
and working out the exterior derivatives
of the $\eta$ we find
\[
d 
\begin{pmatrix}
\eta^i \\
\eta^{\mu} \\
\eta^i_{\mu}
\end{pmatrix}
=
-
\begin{pmatrix}
0 & \eta^i_{\mu} & 0 \\
0 & 0 & 0 \\
0 & 0 & 0  
\end{pmatrix}
\wedge 
\begin{pmatrix}
\eta^i \\
\eta^{\mu} \\
\eta^i_{\mu}
\end{pmatrix}.
\]
This gives
}
\[
d 
\begin{pmatrix}
\omega^i \\
\omega^{\mu} \\
\omega^i_{\mu}
\end{pmatrix}
=
-
\begin{pmatrix}
\omega^i_j & \omega^i_{\nu} + t^i_{\mu \nu} \omega^{\nu} & 0 \\
\omega^{\mu}_j & \omega^{\mu}_{\nu} & 0 \\
\omega^i_{\mu j} & \omega^i_{\mu \nu} & \omega^i_j \delta^{\nu}_{\mu}
- \delta^i_j \omega^{\nu}_{\mu}
\end{pmatrix}
\wedge 
\begin{pmatrix}
\omega^j \\
\omega^{\nu} \\
\omega^j_{\nu}
\end{pmatrix}
\]
where the matrix
\[
\begin{pmatrix}
\omega^i_j & \omega^i_{\nu} + t^i_{\mu \nu} \omega^{\nu} & 0 \\
\omega^{\mu}_j & \omega^{\mu}_{\nu} & 0 \\
\omega^i_{\mu j} & \omega^i_{\mu \nu} & \omega^i_j \delta^{\nu}_{\mu}
- \delta^i_j \omega^{\nu}_{\mu}
\end{pmatrix}
\]
is called the \emph{pseudoconnection 1-form}.
\extraStuff{
It can be described
as follows: first we define the matrix of 1-forms
\[
\Omega = da \, a^{-1}
\]
or, to write it out in full
\[
\begin{pmatrix}
\Omega^i_j & 0 & 0 \\
\Omega^{\mu}_j & \Omega^{\mu}_{\nu} & 0 \\
\Omega^i_{\mu j} & \Omega^i_{\mu \nu} & 
\delta^i_j \Omega^{\nu}_{\mu} - \Omega^i_j \delta^{\nu}_{\mu}
\end{pmatrix}
=
d
\begin{pmatrix}
a^i_k & 0 & 0 \\
a^{\mu}_k & a^{\mu}_{\sigma} & 0 \\
a^i_{\mu k} & a^i_{\mu \sigma} & a^i_k A^{\sigma}_{\mu}  
\end{pmatrix}
\begin{pmatrix}
A^k_j & 0 & 0 \\
A^{\sigma}_j & A^{\sigma}_{\nu} & 0 \\
A^k_{\sigma j} & A^k_{\sigma \nu} & A^k_j a^{\nu}_{\sigma} \end{pmatrix}.
\]
Then we define
\begin{equation} \label{eqn:Prolongation}
\begin{pmatrix}
\omega^i_j \\
\omega^{\mu}_j \\
\omega^{\mu}_{\nu} \\
\omega^i_{\mu j} \\
\omega^i_{\mu \nu}
\end{pmatrix}
=
\begin{pmatrix}
\Omega^i_j \\
\Omega^{\mu}_j \\
\Omega^{\mu}_{\nu} \\
\Omega^i_{\mu j} \\
\Omega^i_{\mu \nu}
\end{pmatrix}
+
\left ( 
\begin{pmatrix}
p^i_{jk} & p^i_{j \sigma} & p^{i \sigma}_{j k} \\
p^{\mu}_{jk} & p^{\mu}_{j \sigma} & p^{\mu \sigma}_{jk} \\
p^{\mu}_{\nu k} & p^{\mu}_{\nu \sigma} & p^{\mu \sigma}_{\nu k} \\ 
p^i_{\mu j k} & p^i_{\mu j \sigma} & p^{i \sigma}_{\mu j k} \\
p^i_{\mu \nu k} & p^i_{\mu \nu \sigma} & p^{i \sigma}_{\mu \nu k}
\end{pmatrix}
+
\begin{pmatrix}
a^i_{jk} & 0 & 0 \\
a^{\mu}_{jk} & a^{\mu}_{j \sigma} & 0 \\
a^{\mu}_{\nu k} & a^{\mu}_{\nu \sigma} & 0 \\
a^i_{\mu j k} & a^i_{\mu j \sigma} & a^i_{kj} \delta^{\sigma}_{\mu}
        - \delta^i_k a^{\sigma}_{\mu j} \\
a^i_{\mu \nu k} & a^i_{\mu \nu \sigma} & - \delta^i_k a^{\sigma}_{\mu \nu}
\end{pmatrix}
\right )
\wedge
\begin{pmatrix}
\omega^k \\
\omega^{\sigma} \\
\omega^k_{\sigma}
\end{pmatrix}
\end{equation}
where the $a^i_{jk}$ etc. satisfy 
the relations
\begin{equation} \label{eqn:aRelations}
\begin{split} 
a^i_{jk} =& a^i_{kj} \\
a^{\mu}_{jk} =& a^{\mu}_{kj} \\
a^{\mu}_{j \sigma} =& a^{\mu}_{\sigma j} \\
a^{\mu}_{\nu \sigma} =& a^{\mu}_{\sigma \nu} \\
a^i_{\mu j k} =& a^i_{\mu k j} \\
a^i_{\mu j \sigma} =& a^i_{\mu \sigma j} \\
a^i_{\mu \nu \sigma} =& a^i_{\mu \sigma \nu}
\end{split}
\end{equation}
and otherwise these new $a$ parameters
are arbitrary; on the other hand the 
$p$ functions are determined entirely
in terms of the previously defined variables
by the equations in table \vref{tbl:prolongation}.
\begin{sidewaystable}
\centering
\begin{xalignat*}{3}
p^i_{jk} =&
\frac{1}{2} a^i_l 
\left (
        A^l_{\epsilon k} A^{\epsilon}_j
        -
        A^l_{\epsilon j} A^{\epsilon}_k
\right ) &
p^i_{j \sigma} =&
-a^i_l \left ( 
        A^l_{\epsilon j} A^{\epsilon}_{\sigma} 
        -       
        A^l_{\epsilon \sigma} A^{\epsilon}_j 
        \right ) &
p^{i \sigma}_{j k} =&
a^{\sigma}_{\epsilon} A^{\epsilon}_j \delta^i_k \\
p^{\mu}_{jk} =& 
        \frac{1}{2} a^{\mu}_l 
        \left (
        A^l_{\epsilon k} A^{\epsilon}_j 
        -       
        A^l_{\epsilon j} A^{\epsilon}_k 
        \right ) &
p^{\mu}_{j \sigma} =&
        - \frac{1}{2} a^{\mu}_l 
        \left (
        A^l_{\epsilon j} A^{\epsilon}_{\sigma} 
        -
        A^l_{\epsilon \sigma} A^{\epsilon}_j 
        \right ) &
p^{\mu \sigma}_{j k} =&
        a^{\mu}_l A^l_k a^{\sigma}_{\epsilon} A^{\epsilon}_j \\
p^{\mu}_{\nu k} =&
        \frac{1}{2} a^{\mu}_l 
        \left ( 
        A^l_{\epsilon k} A^{\epsilon}_{\nu} 
        -       
        A^l_{\epsilon \nu} A^{\epsilon}_k 
        \right ) &
p^{\mu}_{\nu \sigma} =&
        \frac{1}{2} a^{\mu}_j 
        \left (
        A^j_{\epsilon \sigma} A^{\epsilon}_{\nu} 
        -
        A^j_{\epsilon \nu} A^{\epsilon}_{\sigma}
        \right ) &
p^{\mu \sigma}_{\nu k} =&
        a^{\mu}_j A^j_k \delta^{\sigma}_{\nu} \\
p^i_{\mu j k} =& 
a^i_{\mu l}
        \left (
        A^l_{\epsilon k} A^{\epsilon}_j 
        -
        A^l_{\epsilon j} A^{\epsilon}_k
        \right ) &
p^i_{\mu j \sigma} =& 
\frac{1}{2} a^i_{\mu l}
        \left (
        A^l_{\epsilon j} A^{\epsilon}_{\sigma} 
        -
        A^l_{\epsilon \sigma} A^{\epsilon}_l
        \right) &
p^{i \sigma}_{\mu j k} =&
p^i_{jk} \delta^{\sigma}_{\mu}
-
\delta^i_k p^{\sigma}_{j \mu} 
+a^i_{\mu l} A^l_k a^{\sigma}_m A^m_j
\\ 
p^i_{\mu \nu k} =&
- \frac{1}{2} a^i_{\mu j} 
        \left (
        A^j_{\epsilon k} A^{\epsilon}_{\nu} 
        -
        A^j_{\epsilon \nu} A^{\epsilon}_k 
        \right ) &
p^i_{\mu \nu \sigma} =&
\frac{1}{2} a^i_{\mu j} 
        \left (
        A^j_{\epsilon \sigma} A^{\epsilon}_{\nu}
        -
        A^j_{\epsilon \nu} A^{\epsilon}_{\sigma}
        \right ) \\
p^{i \sigma}_{\mu \nu k} =&
-a^i_j \left (
        A^j_{\epsilon k} A^{\epsilon}_{\nu} 
        -
        A^j_{\epsilon \nu} A^{\epsilon}_k 
        \right ) \delta^{\sigma}_{\mu} 
\\ & - \frac{1}{2} \delta^i_k a^{\sigma}_j \left (
        A^j_{\epsilon \nu} A^{\epsilon}_{\mu}
        -
        A^j_{\epsilon \mu} A^{\epsilon}_{\nu}
        \right ) 
\\ & +a^i_{\mu j} A^j_k \delta^{\sigma}_{\nu}
\end{xalignat*}
\caption{Functions defining the prolongation}
\label{tbl:prolongation}
\end{sidewaystable}
Putting this together we find that
\[
t^i_{\mu \nu} =
\frac{1}{2} 
\left (
a^i_{\mu \epsilon} A^{\epsilon}_{\nu}
-
a^i_{\nu \epsilon} A^{\epsilon}_{\mu}
\right )
\]
and that these $t^i_{\mu \nu}$ functions
on $B_0$ vanish precisely on the locus where
\[
a^i_{\mu \epsilon} a^{\mu}_{\nu} =
a^i_{\mu \nu} a^{\mu}_{\epsilon}
\]
which is easily seen to be a smooth
submanifold of $B_0$, say $B$, and
moreover a principal $G$ subbundle of $B_0$,
where $G$ is the subgroup of $G_0$ 
given by matrices $a$ 
satisfying these same equations, i.e.
$G$ is the group of matrices of the form
\[
\begin{pmatrix}
a^i_j & 0 & 0 \\
a^{\mu}_j & a^{\mu}_{\nu} & 0 \\
a^i_{\mu j} & a^i_{\mu \nu} & a^i_j A^{\nu}_{\mu} 
\end{pmatrix}
\text{ with } a^i_{\mu \epsilon} a^{\mu}_{\nu} 
= a^i_{\mu \nu} a^{\mu}_{\epsilon}.
\]
}
On the locus $B=\left(t^i_{\mu \nu}=0\right)$ we find the structure
equations
\begin{equation} \label{eqn:OriginalStruc}
d 
\begin{pmatrix}
\omega^i \\
\omega^{\mu} \\
\omega^i_{\mu}
\end{pmatrix}
=
-
\begin{pmatrix}
\omega^i_j & \omega^i_{\nu} & 0 \\
\omega^{\mu}_j & \omega^{\mu}_{\nu} & 0 \\
\omega^i_{\mu j} & \omega^i_{\mu \nu} & \omega^i_j  \delta^{\nu}_{\mu}
- \delta^i_j \omega^{\nu}_{\mu}
\end{pmatrix}
\wedge 
\begin{pmatrix}
\omega^j \\
\omega^{\nu} \\
\omega^j_{\nu}
\end{pmatrix}
\end{equation}
and
\[
\omega^i_{\mu \nu} = \omega^i_{\nu \mu}.
\]
\extraStuff{
These $\omega^i_j$ etc. (constituting
the connection 1-form)
are now defined up
to addition of the terms $a^i_{jk}$ etc.
which satisfy equations 
\vref{eqn:aRelations} as well as
\[
a^i_{\mu \nu j} = a^i_{\nu \mu j}.
\]
These $a^i_{jk}$ etc. free parameters in the
$\omega^i_j$ etc. 1-forms are the coordinates
of the fibers of the first prolongation $B^{(1)}$.
Under the right $G$ action, the pseudoconnection 
1-form transforms via the adjoint representation:
\[
r_a^*
\begin{pmatrix}
\omega^i_j & \omega^i_{\nu} & 0 \\
\omega^{\mu}_j & \omega^{\mu}_{\nu} & 0 \\
\omega^i_{\mu j} & \omega^i_{\mu \nu} & \omega^i_j  \delta^{\nu}_{\mu}
- \delta^i_j \omega^{\nu}_{\mu}
\end{pmatrix}
=
a^{-1}
\begin{pmatrix}
\omega^i_j & \omega^i_{\nu} & 0 \\
\omega^{\mu}_j & \omega^{\mu}_{\nu} & 0 \\
\omega^i_{\mu j} & \omega^i_{\mu \nu} &  \omega^i_j \delta^{\nu}_{\mu}
- \delta^i_j \omega^{\nu}_{\mu}
\end{pmatrix}
a.
\]
}
\begin{proposition}\label{prop:Diff}
The group of diffeomorphisms of $M$ 
acts transitively on the manifold $B$.
\end{proposition}
\begin{proof}
Given any adapted coframe $U \in B$
we need only show that there is a
system of adapted coordinates $x,y,p$
in which
\[
U = \eta^i \oplus \eta^{\mu} \oplus
\eta^i_{\mu}
\]
at the origin of coordinates.
Let us start by picking any adapted coordinates $x,y,p$.
Certainly we have $U=a \eta$,
for some $a \in G$, no matter what adapted
coordinates we picked. Now change coordinates
by 
\begin{align*}
Y^i &= a^i_j y^j \\
X^{\mu} &= a^{\mu}_j y^j + a^{\mu}_{\nu} x^{\nu} 
\end{align*}
Then in the new coordinates,
\[
\begin{pmatrix}
a^i_j & 0 \\
a^{\mu}_j & a^{\mu}_{\nu} 
\end{pmatrix}
=
\begin{pmatrix}
\delta^i_j & 0 \\
a^{\mu}_j & \delta^{\mu}_{\nu} 
\end{pmatrix}.
\]
We still have the last row to deal with.
Try the change of coordinates
\begin{align*}
Y^i &= y^i \\
X^{\mu} &= x^{\mu} + a^{\mu}_j y^j
\end{align*}
and you find that this accomplishes
the task at hand (only at the origin
of coordinates). Therefore every 
adapted coframe from $B$ arises from
adapted coordinates. Given any two adapted
coframes, take such coordinates near
each of them, and
as diffeomorphism use these coordinate
functions.
\end{proof}
We have therefore discovered the structure equations
(in the sense of Cartan) of the canonical $G$ structure
on the Grassmann bundle.
\extraStuff{
The reader familiar with the theory of $G$ structures
will see that these structure equations look very
like the structure equations of a $\GL{n,\R{}}$ structure
(the equations of a smooth structure),
except that the $\omega^i_{\mu}$ are semibasic,
rather than belonging to the first prolongation.
Indeed the equations of a $\GL{n+d,\R{}}$ structure
on an manifold $M^{n+d}$ are
\begin{equation} \label{eqn:GLNR}
\begin{split}
0 &= d \omega^I + \omega^I_J \wedge \omega^J  \\
0 &= d \omega^I_J + \omega^I_K \wedge \omega^K_J
+ \omega^I_{JK} \wedge \omega^K \\
0 &= d \omega^I_{JK} + \omega^I_L \wedge \omega^L_{JK}
+ \omega^I_{LK} \wedge \omega^L_J 
+ \omega^I_{JL} \wedge \omega^L_K 
+ \omega^I_{JKL} \wedge \omega^L 
\end{split}
\end{equation}
where the capital Roman indices run over first
all $d$ values of the Roman indices, and then all $n$
values of the Greek indices, and these 1-forms
are symmetric in their lower indices. (These equations
are defined on the second prolongation.)
So to produce the equations of the Grassmann bundle,
just split each capital Roman index into two:
a small Roman and a small Greek.
This similarity extends to prolongations of all
orders, so that we can easily use the structure
equations of the prolongations of 
a smooth structure to write down
the structure equations of all prolongations
of the canonical $G$ structure on the
Grassmann bundle. An explanation of this
phenomenon and its generalizations is given
in McKay \cite{McKay:Unpub}. 
For example, it tells us that on the first
prolongation
\begin{align*}
d
\begin{pmatrix}
\omega^i_j & \omega^i_{\nu} \\
\omega^{\mu}_j & \omega^{\mu}_{\nu}
\end{pmatrix}
=&
-
\begin{pmatrix}
\omega^i_k & \omega^i_{\sigma} \\
\omega^{\mu}_k & \omega^{\mu}_{\sigma }
\end{pmatrix}
\wedge
\begin{pmatrix}
\omega^k_j & \omega^k_{\nu} \\
\omega^{\sigma}_j & \omega^{\sigma}_{\nu}
\end{pmatrix}
\\&-
\begin{pmatrix}
\omega^i_{j k} \wedge \omega^k
+ \omega^i_{j \sigma} \wedge \omega^{\sigma} & 
\omega^i_{\nu k} \wedge \omega^k
+ \omega^i_{\nu \sigma} \wedge \omega^{\sigma} \\
\omega^{\mu}_{j k} \wedge \omega^k 
+ \omega^{\mu}_{j \sigma} \wedge \omega^{\sigma} &
\omega^{\mu}_{\nu k} \wedge \omega^k
+
\omega^{\mu}_{\nu \sigma} \wedge \omega^{\sigma}
\end{pmatrix}.
\end{align*}
}
\subsection{Reconstructing the Grassmannian from the
structure equations}
Conversely, suppose that we are given a manifold
$X$ of dimension $n+d+nd$, and on it a coframing
by 1-forms $\xi^i,\xi^{\mu},\xi^i_{\mu}$
and that there exist 1-forms $\xi^i_j$ etc.
so that
\[
d 
\begin{pmatrix}
\xi^i \\
\xi^{\mu} \\
\xi^i_{\mu}
\end{pmatrix}
=
-
\begin{pmatrix}
\xi^i_j & \xi^i_{\nu} & 0 \\
\xi^{\mu}_j & \xi^{\mu}_{\nu} & 0 \\
\xi^i_{\mu j} & \xi^i_{\mu \nu} & \xi^i_j \delta^{\nu}_{\mu}
- \delta^i_j \xi^{\nu}_{\mu}
\end{pmatrix}
\wedge 
\begin{pmatrix}
\xi^j \\
\xi^{\nu} \\
\xi^j_{\nu}
\end{pmatrix}.
\]
Then we have a foliation cut out by the equations
\[
\xi^i=\xi^{\mu}=0.
\]
Suppose that this foliation consists of the
stalks of a submersion,
which is always the case locally. Write this
submersion as $\rho : X \to M$; this will
define the manifold $M$. Let $\Theta$ be 
the plane field on $X$ consisting of tangent
vectors satisfying the equations
\[
\xi^i = 0.
\]
Then we have a map
\[
\phi : X \to \Gro{n}{TM}
\]
defined by 
\[
x \in X \mapsto \phi(x) = \rho'(x) \Theta(x) \in \Gro{n}{TM}.
\]
\begin{proposition} \label{prop:StrucEqns}
The map $\phi : X \to \Gro{n}{TM}$
is a local diffeomorphism, so that under this 
diffeomorphism the 1-forms
$\xi^i,\xi^{\mu},\xi^i_{\mu}$ become
(locally) an adapted coframing of the Grassmann bundle.
\end{proposition}
\ifthenelse{\boolean{abridged}}
  {
The proof is elementary.
  }
  {
\begin{proof} To see that the map $\phi$ is well
defined, note that 
\(
\rho'(x) \Theta(x) 
\)
is an $n$ plane because $\Theta$ is an 
$n+nd$ plane, and the kernel of $\rho'$ is entirely
contained inside $\Theta(x)$, and has dimension $nd$.
We now have established the commutative diagram
\[
\xymatrix{
X \ar[dr]_{\rho} \ar[rr]^{\phi} & & \Gro{n}{TM} \ar[dl]^{\pi} \\
& M &
}
\]

  Pick $\eta^i,\eta^{\mu},\eta^i_{\mu}$ any adapted
coframing on an open subset of $\Gro{n}{TM}$.
The map $\pi : \Gro{n}{TM} \to M$ determines the
plane field $\Theta$ by
\[
\Theta(P) = \pi'(P)^{-1} P.
\]
The same equation holds on $X$:
\[
\Theta(x) = \rho'(x)^{-1}P
\]
where $P = \phi(x)$. This gives
\begin{align*}
\left \{ \phi^* \eta^i = 0 \right \} 
&=
\phi'(x)^{-1} \left \{ \eta^i = 0 \right \} \\
&=
\phi'(x)^{-1} \Theta(P) \\
&=
\phi'(x)^{-1} \pi'(P)^{-1} P \\
&=
\left ( \pi \phi \right )'(x) ^{-1} P \\
&=
\rho'(x)^{-1} P \\
&=
\Theta(x) \\
&= 
\left \{ \xi^i = 0 \right \}
\end{align*}
Therefore there is an invertible matrix $a^i_j$ so that
\[
\phi^* \eta^i = a^i_j \xi^j.
\]
A similar argument establishes that
there is an invertible matrix $a^{\mu}_{\nu}$
so that
\[
\phi^*
\begin{pmatrix}
\eta^i \\
\eta^{\mu}
\end{pmatrix}
=
\begin{pmatrix}
a^i_j & 0 \\
a^{\mu}_j & a^{\mu}_{\nu} 
\end{pmatrix}
\begin{pmatrix}
\xi^j \\
\xi^{\nu}
\end{pmatrix}.
\]
Taking the exterior derivatives of these
equations, and plugging in the structure
equations, we find by Cartan's lemma that there are functions
$a^i_{\mu j}, a^i_{\mu \nu}$ so that
\[
\phi^* \eta^i_{\mu}
=
a^i_{\mu j} \xi^j + a^i_{\mu \nu} \xi^{\nu} + a^i_j \xi^j_{\nu} A^{\nu}_{\mu}
\]
or in other words that
\[
\phi^* 
\begin{pmatrix}
\eta^i \\
\eta^{\mu} \\
\eta^i_{\mu}
\end{pmatrix}
=
\begin{pmatrix}
a^i_j & 0 & 0 \\
a^{\mu}_j & a^{\mu}_{\nu} & 0 \\
a^i_{\mu j} & a^i_{\mu \nu} & a^i_j A^{\nu}_{\mu}
\end{pmatrix}
\begin{pmatrix}
\xi^j \\
\xi^{\nu} \\
\xi^j_{\nu}
\end{pmatrix}
\]
i.e. a change of adapted coframing. Therefore
$\phi'$ has full rank, and the result is clear.
\end{proof}
}

\subsection{Complex notation}
We wish to describe the same structure
equations using a complex notation, 
assuming that the manifold $M$ is of
even dimension (say $2(n+d)$), and that the planes from
which the Grassmann bundle is composed
are also of even dimension (say $2n$).
Then we can pick any complex structure on the
vector space $V$, so that the chosen subspace $P_0$ 
is a complex subspace. Using a complex basis
of $V$ instead of a real one, we find that the same
structure equations hold that we already had,
but we have only to split the 1-forms, now complex
valued, into complex linear and antilinear parts
on the complexified tangent bundle. 
\extraStuff{
More precisely,
we find that we can write our structure equations
as in table~\vref{tbl:strucEqnsOne}, where
\begin{xalignat*}{3}
\omega^i_{\mu \nu} =& \omega^i_{\nu \mu} &
\omega^{\bari}_{\mu \nu} =& \omega^{\bari}_{\nu \mu} &
\omega^i_{\mu \bar{\nu}} =& \omega^i_{\bar{\nu} \mu} \\
\omega^{\bari}_{\mu \bar{\nu}} =& \omega^{\bari}_{\bar{\nu} \mu} &
\omega^i_{\bar{\mu} \bar{\nu}} =& \omega^i_{\bar{\nu} \bar{\mu}} &
\omega^{\bari}_{\bar{\mu} \bar{\nu}} =& 
        \omega^{\bari}_{\bar{\nu} \bar{\mu}}.
\end{xalignat*}
Now instead of splitting the capital Roman
indices from equation \vref{eqn:GLNR} into
small Roman and small Greek, we have to
split into small Roman, small Roman barred,
small Greek and small Greek barred.
% \begin{landscape}
% \newpage
\begin{sidewaystable}
\centering
\[
d 
\begin{pmatrix}
\omega^i \\
\omega^{\bar i} \\
\omega^{\mu} \\
\omega^{\bar \mu} \\
\omega^i_{\mu} \\
\omega^{\bar i}_{\bar \mu} \\
\omega^i_{\bar \mu} \\
\omega^{\bar i}_{\mu}
\end{pmatrix}
=
-
\begin{pmatrix}
\omega^i_j & \omega^i_{\bar j} & \omega^i_{\nu} & \omega^i_{\bar \nu} 
& 0 & 0 & 0 & 0 \\
\omega^{\bar i}_j & \omega^{\bar i}_{\bar j} & 
\omega^{\bar i}_{\nu} & \omega^{\bar i}_{\bar \nu} 
& 0 & 0 & 0 & 0 \\
\omega^{\mu}_j & 
\omega^{\mu}_{\bar j} &
\omega^{\mu}_{\nu} &
\omega^{\mu}_{\bar \nu} &
0 & 0 & 0 & 0 \\
\omega^{\bar \mu}_j & 
\omega^{\bar \mu}_{\bar j} &
\omega^{\bar \mu}_{\nu} &
\omega^{\bar \mu}_{\bar \nu} &
0 & 0 & 0 & 0 \\
\omega^i_{\mu j} & 
\omega^i_{\mu \barj} &
\omega^i_{\mu \nu} &
\omega^i_{\mu \bar{\nu}} &
\omega^i_j \delta^{\nu}_{\mu} 
- \delta^i_j \omega^{\nu}_{\mu} &
0 & 
-\delta^i_j \omega^{\bar \nu}_{\mu} & 
\omega^i_{\barj} \delta^{\nu}_{\mu}  \\
\omega^{\bar i}_{\bar{\mu} j} &
\omega^{\bar i}_{\bar{\mu} \barj} &
\omega^{\bar i}_{\bar{\mu} \nu} &
\omega^{\bar i}_{\bar{\mu} \bar{\nu}} &
0 &
\omega^{\bari}_{\barj} \delta^{\bar{\nu}}_{\bar{\mu}} 
-
\delta^{\bari}_{\barj} \omega^{\bar{\nu}}_{\bar{\mu}}
&
\omega^{\bari}_{j} \delta^{\bar{\nu}}_{\bar{\mu}}  &
-\delta^{\bari}_{\barj} \omega^{\nu}_{\bar{\mu}} \\
\omega^{i}_{\bar{\mu} j} &
\omega^{i}_{\bar{\mu} \barj} &
\omega^{i}_{\bar{\mu} \nu} &
\omega^{i}_{\bar{\mu} \bar{\nu}} &
- \delta^i_j \omega^{\nu}_{\bar{\mu}} &
\omega^i_{\barj} \delta^{\bar{\nu}}_{\bar{\mu}} &
\omega^i_j \delta^{\bar{\nu}}_{\bar{\mu}} 
-
\delta^i_j \omega^{\bar{\nu}}_{\bar{\mu}} &
0 \\
\omega^{\bari}_{\mu j} &
\omega^{\bari}_{\mu \barj} &
\omega^{\bari}_{\mu \nu} &
\omega^{\bari}_{\mu \bar{\nu}} &
\omega^{\bari}_j \delta^{\nu}_{\mu} &
- \delta^{\bari}_{\barj} \omega^{\bar{\nu}}_{\mu} & 
0 & 
\omega^{\bari}_{\barj} \delta^{\nu}_{\mu} 
-
\delta^{\bari}_{\barj} \omega^{\nu}_{\mu}
\end{pmatrix}
\wedge 
\begin{pmatrix}
\omega^j \\
\omega^{\bar j} \\
\omega^{\nu} \\
\omega^{\bar{\nu}} \\
\omega^j_{\nu} \\
\omega^{\bar j}_{\bar \nu} \\
\omega^j_{\bar \nu} \\
\omega^{\bar j}_{\nu}
\end{pmatrix}
\]
\caption{The structure equations on $\Gro{2n}{TM}$ in complex notation}
\label{tbl:strucEqnsOne}
\end{sidewaystable}
% \end{landscape}
}
Take any system of complex valued coordinates
\[
z^{\mu}, w^i
\]
on $M$ near a point $m \in M$.\extraStuff{\footnote{These coordinates are 
not in any sense required to be holomorphic; there is no
sense in which they could be, since the manifold $M$ is
only a real manifold of even dimension.}}
Every $2n$-plane on which the $dz^{\mu}$ and $dz^{\bar{\mu}}$
are linearly independent is described by an
equation like
\[
dw^i = p^i_{\mu} \, dz^{\mu} + p^i_{\bar{\mu}} \, dz^{\bar{\mu}}
\]
%% \begin{align*}
%% dw^i =& p^i_{\mu} \, dz^{\mu} + p^i_{\bar{\mu}} \, dz^{\bar{\mu}} \\
%% dw^{\bari} =& p^{\bari}_{\mu} \, dz^{\mu} + p^{\bari}_{\bar{\mu}} \,
%% dz^{\bar{\mu}}
%% \end{align*}
\extraStuff{
with the convention that
\begin{align*}
p^{\bari}_{\bar{\mu}} =& \overline{p^i_{\mu}} \\
p^{\bari}_{\mu} =& \overline{p^i_{\bar{\mu}}}
\end{align*}
are complex conjugates. 
}
Thus the numbers
\[
z^{\mu}, w^i, p^i_{\mu}, p^i_{\bar{\mu}}
\]
provide complex valued coordinates on $\Gro{2k}{TM}$
near the $2n$-plane $dw=0$.
\extraStuff{We can describe 
a section of the bundle $B$ as
\begin{align*}
\eta^i =& dw^i - p^i_{\mu} \, dz^{\mu}
        -p^i_{\bar{\mu}} \, dz^{\bar{\mu}} \\
%\eta^{\bari} =& dw^{\bari} - p^{\bari}_{\mu} \, dz^{\mu}
%        -p^{\bari}_{\bar{\mu}} \, dz^{\bar{\mu}} \\
\eta^{\mu} =&
dz^{\mu} \\
%\eta^{\bar{\mu}} =&
%dz^{\bar{\mu}} \\
\eta^i_{\mu} =&
dp^{i}_{\mu} \\
%\eta^{\bari}_{\bar{\mu}} =&
%dp^{\bari}_{\bar{\mu}} \\
\eta^i_{\bar{\mu}} =&
dp^i_{\bar{\mu}} \\
%\eta^{\bari}_{\mu} =&
%dp^{\bari}_{\mu}
\end{align*}
The soldering 1-form is written
in table~\vref{tbl:SolderingOne}.
% \begin{landscape}
\begin{table}\label{tbl:SolderingOne}
\[
\begin{pmatrix}
\omega^i \\
\omega^{\bari} \\
\omega^{\mu} \\
\omega^{\bar{\mu}} \\
\omega^i_{\mu} \\
\omega^{\bari}_{\bar{\mu}} \\
\omega^i_{\bar{\mu}} \\
\omega^{\bari}_{\mu}
\end{pmatrix}
=
\begin{pmatrix}
a^i_j & a^i_{\barj} & 0 & 0 & 0 & 0 & 0 & 0 \\
a^{\bari}_j & a^{\bari}_{\barj} & 0 & 0 & 0 & 0 & 0 & 0 \\
a^{\mu}_j & a^{\mu}_{\barj} & a^{\mu}_{\nu} & a^{\mu}_{\bar{\nu}} 
& 0 & 0 & 0 & 0 \\
a^{\bar{\mu}}_j & a^{\bar{\mu}}_{\barj} & 
a^{\bar{\mu}}_{\nu} & a^{\bar{\mu}}_{\bar{\nu}} & 
0 & 0 & 0 & 0 
\\
a^{i}_{\mu j} & a^i_{\mu \barj} &
a^{i}_{\mu \nu} & a^i_{\mu \bar{\nu}} &
a^{i}_{j} A^{\nu}_{\mu} & a^i_{\barj} A^{\bar{\nu}}_{\mu} & 
a^{i}_{j} A^{\bar{\nu}}_{\mu} & a^i_{\barj} A^{\nu}_{\mu}  
\\
a^{\bari}_{\bar{\mu} j} & a^{\bari}_{\bar{\mu} \barj} &
a^{\bari}_{\bar{\mu} \nu} & a^{\bari}_{\bar{\mu} \bar{\nu}} &
a^{\bari}_{j} A^{\nu}_{\bar{\mu}} & 
a^{\bari}_{\barj} A^{\bar{\nu}}_{\bar{\mu}} &
a^{\bari}_{j} A^{\bar{\nu}}_{\bar{\mu}} & 
a^{\bari}_{\barj} A^{\nu}_{\bar{\mu}}
\\
a^{i}_{\bar{\mu} j} & a^{i}_{\bar{\mu} \barj} &
a^i_{\bar{\mu} \nu} & a^i_{\bar{\mu} \bar{\nu}} &
a^{i}_{j} A^{\nu}_{\bar{\mu}} & 
a^{i}_{\barj} A^{\bar{\nu}}_{\bar{\mu}} &
a^{i}_{j} A^{\bar{\nu}}_{\bar{\mu}} & 
a^{i}_{\barj} A^{\nu}_{\bar{\mu}}
\\
a^{\bari}_{\mu j} & a^{\bari}_{\mu \barj} &
a^{\bari}_{\mu \nu} & a^{\bari}_{\mu \bar{\nu}} &
a^{\bari}_{j} A^{\nu}_{\mu} & 
a^{\bari}_{\barj} A^{\bar{\nu}}_{\mu} &
a^{\bari}_{j} A^{\bar{\nu}}_{\mu} & 
a^{\bari}_{\barj} A^{\nu}_{\mu}
\end{pmatrix}
\begin{pmatrix}
\eta^j \\
\eta^{\barj} \\
\eta^{\nu} \\
\eta^{\bar{\nu}} \\
\eta^i_{\nu} \\
\eta^{\barj}_{\bar{\nu}} \\
\eta^i_{\bar{\nu}} \\
\eta^{\barj}_{\nu}
\end{pmatrix}
\]
\caption{Soldering 1-form on $\Gro{2n}{TM}$ in
complex notation}
\end{table}
% \end{landscape}
}
\section{Complex geometry}
Given a complex structure on $M$, we can look at the
complex Grassmann bundle
\[
\iota : \GrC{k}{TM} \subset \Gro{2k}{TM}
\]
which has structure equations given by
the same reasoning in purely holomorphic
terms, so
\[
d 
\begin{pmatrix}
\omega^i \\
\omega^{\mu} \\
\omega^i_{\mu}
\end{pmatrix}
=
-
\begin{pmatrix}
\omega^i_j & \omega^i_{\mu} & 0 \\
\omega^{\mu}_j & \omega^{\mu}_{\nu} & 0 \\
\omega^i_{\mu j} & \omega^i_{\mu \nu} &  \omega^i_j \delta^{\nu}_{\mu}
- \delta^i_j \omega^{\nu}_{\mu}
\end{pmatrix}
\wedge 
\begin{pmatrix}
\omega^i \\
\omega^{\mu} \\
\omega^i_{\mu}
\end{pmatrix}
\]
and the conjugates of these equations. These
hold on a bundle 
\[
\Pi_{\C{}} : B_{\C{}} \to \GrC{k}{TM}
\]
constructed by carrying out the same process
as before, but using only complex linear data.
As in proposition 
\vref{prop:Diff}, the local biholomorphisms
of $M$ act transitively on $B_{\C{}}$. By
the same argument as in proposition \vref{prop:StrucEqns},
these structure equations determine the local
geometry of the holomorphic Grassmann bundle.
The map $\iota$ into the real Grassmann bundle
allows us to pullback the bundle $B$ to
$\GrC{k}{TM}$, and also to map $B_{\C{}}$
into the pullback bundle.
\[
\xymatrix{
B_{\C{}} \ar[dr]_{\Pi_{\C{}}} \ar[r] 
& \iota^* B \ar[r] \ar[d] & B \ar[d]^{\Pi} \\
                            & \GrC{k}{TM} \ar[r]       & \Gro{2k}{TM} \\ 
}
\]
Pulling back the 1-forms from $B$ we find 
$\omega^i$ pulls back to $\omega^i$, etc.,
i.e. that 
the 1-forms on $B$ 
become the holomorphic and conjugate
holomorphic 1-forms on $B^{\C{}}$, except
for those 1-forms which have mixed indices,
i.e. both barred and unbarred indices. These
all vanish. For example, on $B_{\C{}}$
\(
\omega^i_{\bar{\mu}} %=\omega^{\bari}_{\mu} 
= 0.
\)

\section{Analogues of complex geometry}

A differential equation of the type we are
studying imposes itself in this picture as
a submanifold $E \subset \Gro{2n}{TM}$,
via the equations
\[
p^i_{\bar{\mu}} = F^i_{\bar{\mu}}.
\]
Globally, let us assume only that we have an
immersed submanifold
\[
\phi : E^{2(n+d+nd)} \to \Gro{2n}{TM}.
\]
%Dimensions:
%\[
%\xymatrix{
%E^{2(n+d+nd)} \ar[rr]^{\phi} \ar[dr] & & \Gro{2n}{TM}^{2(n+d)+4nd} \ar[dl] \\
% & M^{2(n+d)}
%}
%\]
and that the composition mapping $E \to M$
is a submersion, and that near any point there 
are local coordinates in which $E$ is represented
by equations on the $p^i_{\bar{\mu}}$ (which is
just a nondegeneracy condition).
The manifold $E$ has the same real dimension as the complex
Grassmann bundle, and plays an analogous
role.
We can pull the bundle $B$ back to $E$
via $\phi$.
\extraStuff{
\[
\xymatrix{
\phi^*B \ar[r] \ar[d] & B \ar[d] \\
E \ar[r] & \Gro{2n}{TM}
}
\]
where $\phi^*B$ is (the definition of pullback)
the set of pairs $(e,P)$ so that $e \in E$
and $P \in B$ belongs to the fiber of $B$
over $\phi(e)$. 
}
We find however that on this 
principal $G$ bundle $\phi^* B \to E$ 
the soldering 1-forms 
\[
\omega^i, \omega^{\bari}, \omega^{\mu}, \omega^{\bar{\mu}},
\omega^i_{\mu}, \omega^{\bari}_{\bar{\mu}},
\omega^i_{\bar{\mu}}, \omega^{\bari}_{\mu}
\]
can no longer be independent, because they
are semibasic

Looking at adapted coordinates, we find that the
1-forms $\omega^i,\omega^{\bari},\omega^{\mu},\omega^{\bar{\mu}}$
are semibasic for the projection to $M$, and 
they must remain independent on $E$.
In the adapted coframing $\eta$
we find relations
\begin{align*}
\eta^i_{\bar{\mu}} =&
\left ( 
\pd{p^i_{\bar{\mu}}}{z^{\nu}} 
 +
 \pd{p^i_{\bar{\mu}}}{w^j} p^j_{\nu}
 \right )
 \eta^{\nu}
%  \\& 
 +
 \left ( 
 \pd{p^i_{\bar{\mu}}}{z^{\bar{\nu}}} 
 +
 \pd{p^i_{\bar{\mu}}}{w^j} p^j_{\bar{\nu}}
 \right )
 \eta^{\bar{\nu}}
 \\
 &+
 \pd{p^i_{\bar{\mu}}}{w^j} \eta^j
 +
 \pd{p^i_{\bar{\mu}}}{w^{\barj}} \eta^{\barj}
%  \\&
 +
 \pd{p^i_{\bar{\mu}}}{p^j_{\nu}} \eta^j_{\nu}
 +
 \pd{p^i_{\bar{\mu}}}{p^{\barj}_{\bar{\nu}}} \eta^{\barj}_{\bar{\nu}}
\end{align*}
%% \begin{align*}
%% \eta^i_{\bar{\mu}} =&
%% \left ( 
%% \pd{p^i_{\bar{\mu}}}{z^{\nu}} 
%% +
%% \pd{p^i_{\bar{\mu}}}{w^j} p^j_{\nu}
%% \right )
%% \eta^{\nu}
%% % \\& 
%% +
%% \left ( 
%% \pd{p^i_{\bar{\mu}}}{z^{\bar{\nu}}} 
%% +
%% \pd{p^i_{\bar{\mu}}}{w^j} p^j_{\bar{\nu}}
%% \right )
%% \eta^{\bar{\nu}}
%% \\
%% &+
%% \pd{p^i_{\bar{\mu}}}{w^j} \eta^j
%% +
%% \pd{p^i_{\bar{\mu}}}{w^{\barj}} \eta^{\barj}
%% % \\&
%% +
%% \pd{p^i_{\bar{\mu}}}{p^j_{\nu}} \eta^j_{\nu}
%% +
%% \pd{p^i_{\bar{\mu}}}{p^{\barj}_{\bar{\nu}}} \eta^{\barj}_{\bar{\nu}}
%% \\
%% \eta^{\bari}_{\mu} =&
%% \left ( 
%% \pd{p^{\bari}_{\mu}}{z^{\nu}} 
%% +
%% \pd{p^{\bari}_{\mu}}{w^j} p^j_{\nu}
%% \right )
%% \eta^{\nu}
%% % \\&
%% +
%% \left ( 
%% \pd{p^{\bari}_{\mu}}{z^{\bar{\nu}}} 
%% +
%% \pd{p^{\bari}_{\mu}}{w^j} p^j_{\bar{\nu}}
%% \right )
%% \eta^{\bar{\nu}}
%% \\ &+
%% \pd{p^{\bari}_{\mu}}{w^j} \eta^j
%% +
%% \pd{p^{\bari}_{\mu}}{w^{\barj}} \eta^{\barj}
%% % \\&
%% +
%% \pd{p^{\bari}_{\mu}}{p^j_{\nu}} \eta^j_{\nu}
%% +
%% \pd{p^{\bari}_{\mu}}{p^{\barj}_{\bar{\nu}}} \eta^{\barj}_{\bar{\nu}}
%% \end{align*}
So for coframes from $\phi^* B$ which
are close enough to this adapted coframing,
we find that we can solve for the ${}^i_{\bar{\mu}}$
and ${}^{\bari}_{\mu}$ entries of the coframe
in terms of the other entries of the coframe.
This holds on a dense open subset of $\phi^* B$
(and on a Zariski open subset of each fiber).
Therefore the soldering 1-forms on that dense
open subset of $\phi^* B$ satisfy equations
\[
\omega^i_{\bar{\mu}} =
t^i_{\bar{\mu} j} \omega^j 
+
t^i_{\bar{\mu} \barj} \omega^{\barj} 
+
t^i_{\bar{\mu} \nu} \omega^{\nu} 
+
t^i_{\bar{\mu} \bar{\nu}} \omega^{\bar{\nu}} 
+
t^{i \nu}_{\bar{\mu} j} \omega^j_{\nu} 
+
t^{i \bar{\nu}}_{\bar{\mu} \barj} \omega^{\barj}_{\bar{\nu}}  \\
\]
%% \begin{align*}
%% \omega^i_{\bar{\mu}} =&
%% t^i_{\bar{\mu} j} \omega^j 
%% +
%% t^i_{\bar{\mu} \barj} \omega^{\barj} 
%% +
%% t^i_{\bar{\mu} \nu} \omega^{\nu} 
%% +
%% t^i_{\bar{\mu} \bar{\nu}} \omega^{\bar{\nu}} 
%% +
%% t^{i \nu}_{\bar{\mu} j} \omega^j_{\nu} 
%% +
%% t^{i \bar{\nu}}_{\bar{\mu} \barj} \omega^{\barj}_{\bar{\nu}}  \\
%% \omega^{\bari}_{\mu} =&
%% t^{\bari}_{\mu j} \omega^j 
%% +
%% t^{\bari}_{\mu \barj} \omega^{\barj} 
%% +
%% t^{\bari}_{\mu \nu} \omega^{\nu} 
%% +
%% t^{\bari}_{\mu \bar{\nu}} \omega^{\bar{\nu}} 
%% +
%% t^{\bari \nu}_{\mu j} \omega^j_{\nu} 
%% +
%% t^{\bari \bar{\nu}}_{\bar{\mu} \barj} \omega^{\barj}_{\bar{\nu}}.
%% \end{align*}
Using the equation
\[
r_g^* \omega = g^{-1} \omega
\]
for the right action of the group $G$ on
the bundle $\phi^* B$, we find that we can
arrange the equations
\[
0 = 
t^i_{\bar{\mu} j} =
t^i_{\bar{\mu} \barj} =
t^i_{\bar{\mu} \nu} =
t^i_{\bar{\mu} \bar{\nu}} + t^i_{\bar{\nu} \bar{\mu}} 
= t^{i \nu}_{\bar{\mu} i} 
= t^{i \bar{\mu}}_{\bar{\mu} \barj} 
= t^{\bari \mu}_{\mu j} 
= t^{\bari \bar{\nu}}_{\mu \bari}.
\]
%% \begin{align*}
%% t^i_{\bar{\mu} j} =
%% t^i_{\bar{\mu} \barj} =
%% t^i_{\bar{\mu} \nu} =
%% t^i_{\bar{\mu} \bar{\nu}} + t^i_{\bar{\nu} \bar{\mu}} 
%% = 0 \\
%% t^{\bari}_{\mu j} =
%% t^{\bari}_{\mu \barj} =
%% t^{\bari}_{\mu \bar{\nu}} =
%% t^{\bari}_{\mu \nu} + t^{\bari}_{\nu \mu} 
%% = 0 \\
%% \end{align*}
%% and 
%% \begin{align*}
%% t^{i \nu}_{\bar{\mu} i} &= 0 \\
%% t^{i \bar{\mu}}_{\bar{\mu} \barj} &= 0 \\
%% t^{\bari \mu}_{\mu j} &= 0 \\ 
%% t^{\bari \bar{\nu}}_{\mu \bari} &= 0.
%% \end{align*}
Thus 
\begin{equation} \label{eqn:Torsion}
%\begin{split}
\omega^i_{\bar{\mu}} =
t^i_{\bar{\mu} \bar{\nu}} \omega^{\bar{\nu}} 
+
t^{i \nu}_{\bar{\mu} j} \omega^j_{\nu} 
+
t^{i \bar{\nu}}_{\bar{\mu} \barj} \omega^{\barj}_{\bar{\nu}}  
%\\
%\omega^{\bari}_{\mu} =&
%t^{\bari}_{\mu \nu} \omega^{\nu} 
%+
%t^{\bari \nu}_{\mu j} \omega^j_{\nu} 
%+
%t^{\bari \bar{\nu}}_{\mu \barj} \omega^{\barj}_{\bar{\nu}}.
%\end{split}
\end{equation}
The subset of $\phi^* B$ on which these equations are
satisfied, call it $B_1$, is a principal $G_1$ 
subbundle, where $G_1$ is a certain subgroup of $G$.

  The differential equations \vref{eqn:Eqn}
can now be written in terms of any adapted
coframing $\eta$ as 
\[
\eta^i = 0 
\]
(and conjugate)
which, when differentiated gives the tableau
\[
d \eta^i = - \eta^i_{\mu} \wedge \eta^{\mu}
- 
\left ( t^i_{\bar{\mu} \bar{\nu}} \eta^{\bar{\nu}} 
+ t^{i \nu}_{\bar{\mu} j} \eta^j_{\nu}
+ t^{i \bar{\nu}}_{\bar{\mu} \barj} 
\eta^{\barj}_{\bar{\nu}}
\right )
\wedge 
\eta^{\bar{\mu}}
\]
where the $t$ terms are pulled back
from the bundle $B_1$. We see that 
a Cauchy--Riemann tableau can emerge 
only if we find a way to eliminate
these $t$ terms, by change of coframing.
On the other hand, working out how these
terms transform under the structure
group $G$, 
it is easy to see that
if they don't vanish at a point, then
they don't vanish anywhere on the fiber
of $B$ through that point.

\begin{lemma} The differential
equations \vref{eqn:Eqn} have Cauchy--Riemann
tableau precisely when
\[
0 =  t^i_{\bar{\mu} \bar{\nu}}
= t^{i \nu}_{\bar{\mu} j}
= t^{i \bar{\nu}}_{\bar{\mu} \barj}.
\]
\end{lemma}
Note that if there is only one $z$ variable
and one $w$ variable in equation \ref{eqn:Eqn}
then all of these equations are automatically
satisfied, because of the relations in equation
\vref{eqn:Torsion}.

Henceforth we will assume that our 
differential equations have Cauchy--Riemann
tableau. 
If the functions
$F$ are real analytic, then this implies
(by the Cartan--K{\"a}hler theorem) that there
are local solutions $w(z,\bar{z})$ to 
\ref{eqn:Eqn}
with the same degree of
generality as the Cauchy--Riemann equations.
More precisely,
\begin{proposition} Suppose that $E \subset \Gro{2n}{TM}$
is a real analytic immersed submanifold
determining a system of differential
equations with Cauchy--Riemann tableau.
Given a real $n$-plane, $P \subset T_m M$,
call it \emph{$E$-admissible} if it lies inside
a $2n$-plane $\hat{P} \subset T_m M$
belonging to the submanifold $E \subset \Gro{2n}{TM}$.
There is a discrete set of such $2n$-planes 
$\hat{P}$. Call them the integral
extensions of $P$. The admissible $n$-planes form
an open subset of $\Gro{n}{TM}$.
Every real analytic immersed 
submanifold $\Sigma$ of $M$
of dimension $n$ whose tangent spaces are
admissible lies in a immersed integral
manifold of $E$, i.e. in a submanifold
$\hat{\Sigma} \subset M$ whose tangent
spaces belong to $E$. The largest such
integral manifold $\hat{\Sigma}$ is
uniquely determined by choice of $\Sigma$
(or the infinite jet of $\Sigma$ at one point)
and choice of one integral extension
$\hat{P}$ of one tangent space $P=T_x \Sigma$
so that $\hat{P} = T_x \hat{\Sigma}$.
\end{proposition}
Note that if the fibers of $E \to M$
are compact, then every $n$-plane is
admissible, a well posed infinitesimal
Cauchy problem. In general, noncompactness
of the fibers will lead to integral 
manifolds ``running off the edge'',
i.e. inextendable to any larger integral
manifold, even at smooth points.
\begin{proof} This is an immediate
consequence of the Cartan--K{\"a}hler
theorem; see Cartan \cite{Cartan:1945}.
\end{proof}

\extraStuff{
Now it is easy to calculate that
the subgroup $G_1$ of $G$ leaving the
subbundle $B_1$ on which $\omega^i_{\bar{\mu}}
= \omega^{\bari}_{\mu}=0$ invariant
is precisely the subgroup
consisting of those matrices with 
\[
a^i_{\bar{\mu} j} = a^{\bari}_{\mu j} =
a^i_{\bar{\mu} \barj} = a^{\bari}_{\mu \barj} =
a^i_{\bar{\mu} \nu} = a^{\bari}_{\mu \nu} =
a^i_{\bar{\mu} \bar{\nu}} = a^{\bari}_{\mu \bar{\nu}} =
a^{\nu}_{\bar{\mu}} = a^{\bari}_j = a^i_{\barj} =
a^{\bar{\nu}}_{\mu} = 0
\]
i.e. the matrices of the form
\[
\begin{pmatrix}
a^i_j & 0 & 0 & 0 & 0 & 0 & 0 & 0 \\
0 & a^{\bari}_{\barj} & 0 & 0 & 0 & 0 & 0 & 0 \\
a^{\mu}_j & a^{\mu}_{\barj} & a^{\mu}_{\nu} & 0 
& 0 & 0 & 0 & 0 \\
a^{\bar{\mu}}_j & a^{\bar{\mu}}_{\barj} & 
0 & a^{\bar{\mu}}_{\bar{\nu}} & 
0 & 0 & 0 & 0 
\\
a^{i}_{\mu j} & a^i_{\mu \barj} &
a^{i}_{\mu \nu} & 0 &
a^{i}_{j} A^{\nu}_{\mu} & 0 & 
0 & 0  
\\
a^{\bari}_{\bar{\mu} j} & a^{\bari}_{\bar{\mu} \barj} &
0 & a^{\bari}_{\bar{\mu} \bar{\nu}} &
0 & 
a^{\bari}_{\barj} A^{\bar{\nu}}_{\bar{\mu}} &
0 & 
0
\\
0 & 0 &
0 & 0 &
0 & 
0 &
a^{i}_{j} A^{\bar{\nu}}_{\bar{\mu}} & 
0
\\
0 & 0 &
0 & 0 &
0 & 
0 &
0 & 
a^{\bari}_{\barj} A^{\nu}_{\mu}
\end{pmatrix}
\]

We may view the results so far as saying
that $B_1$ is an integral
manifold of the differential system 
\[
\omega^i_{\bar{\mu}}=\omega^{\bari}_{\mu}=0
\]
on the manifold $B$. As a consequence, 
\begin{equation} \label{eqn:POPO}
\begin{split}
0 =& d \omega^i_{\bar{\mu}} \\
  =& - \omega^i_{\bar{\mu} j} \wedge \omega^j 
     - \omega^i_{\bar{\mu} \barj} \wedge \omega^{\barj}
     - \omega^i_{\bar{\mu} \nu} \wedge \omega^{\nu}
     - \omega^i_{\bar{\mu} \bar{\nu}} \wedge \omega^{\bar{\nu}} \\
     &+ \delta^i_j \omega^{\nu}_{\bar{\mu}} \wedge \omega^j_{\nu}
     - \omega^i_{\barj} \delta^{\bar{\nu}}_{\bar{\mu}} \wedge 
                \omega^{\barj}_{\bar{\nu}}
\end{split}
\end{equation}
and the conjugates of these equations. By Cartan's
lemma:
\begin{equation} \label{eqn:BIG}
\begin{pmatrix}
\omega^i_{\bar{\mu} j} \\
\omega^i_{\bar{\mu} \barj} \\
\omega^i_{\bar{\mu} \nu} \\
\omega^i_{\bar{\mu} \bar{\nu}} \\
- \delta^i_j \omega^{\nu}_{\bar{\mu}} \\
\omega^i_{\barj} \delta^{\bar{\nu}}_{\bar{\mu}} 
\end{pmatrix}
=
\begin{pmatrix}
t^i_{\bar{\mu} jk} &
t^i_{\bar{\mu} j \bar{k}} &
t^i_{\bar{\mu} j \sigma} &
t^i_{\bar{\mu} j \bar{\sigma}} &
t^{i \sigma}_{\bar{\mu} jk} &
t^{i \bar{\sigma}}_{\bar{\mu} j \bar{k}} \\
t^i_{\bar{\mu} \barj k} &
t^i_{\bar{\mu} \barj \bar{k}} &
t^i_{\bar{\mu} \barj \sigma} &
t^i_{\bar{\mu} \barj \bar{\sigma}} &
t^{i \sigma}_{\bar{\mu} \barj k} &
t^{i \bar{\sigma}}_{\bar{\mu} \barj \bar{k}} \\
t^i_{\bar{\mu} \nu k} &
t^i_{\bar{\mu} \nu \bar{k}} &
t^i_{\bar{\mu} \nu \sigma} &
t^i_{\bar{\mu} \nu \bar{\sigma}} &
t^{i \sigma}_{\bar{\mu} \nu k} &
t^{i \bar{\sigma}}_{\bar{\mu} \nu \bar{k}} \\
t^i_{\bar{\mu} \bar{\nu} k} &
t^i_{\bar{\mu} \bar{\nu} \bar{k}} &
t^i_{\bar{\mu} \bar{\nu} \sigma} &
t^i_{\bar{\mu} \bar{\nu} \bar{\sigma}} &
t^{i \sigma}_{\bar{\mu} \bar{\nu} k} &
t^{i \bar{\sigma}}_{\bar{\mu} \bar{\nu} \bar{k}} \\
t^{i \nu}_{\bar{\mu} jk} &
t^{i \nu}_{\bar{\mu} j \bar{k}} &
t^{i \nu}_{\bar{\mu} j \sigma} &
t^{i \nu}_{\bar{\mu} j \bar{\sigma}} &
t^{i \nu \sigma}_{\bar{\mu} jk} &
t^{i \nu \bar{\sigma}}_{\bar{\mu} j \bar{k}} \\
t^{i \bar{\nu}}_{\bar{\mu} \barj k} &
t^{i \bar{\nu}}_{\bar{\mu} \barj \bar{k}} &
t^{i \bar{\nu}}_{\bar{\mu} \barj \sigma} &
t^{i \bar{\nu}}_{\bar{\mu} \barj \bar{\sigma}} &
t^{i \bar{\nu} \sigma}_{\bar{\mu} \barj k} &
t^{i \bar{\nu} \bar{\sigma}}_{\bar{\mu} \barj \bar{k}} 
\end{pmatrix}
\begin{pmatrix}
\omega^k \\
\omega^{\bar{k}} \\
\omega^{\sigma} \\
\omega^{\bar{\sigma}} \\
\omega^k_{\sigma} \\
\omega^{\bar{k}}_{\bar{\sigma}}
\end{pmatrix}
\end{equation}
and conjugates, with the $t$ functions being
symmetric in their indices so 
that this is a symmetric matrix.
The $\delta$ factors above show that
\[
\begin{pmatrix}
\omega^{\nu}_{\bar{\mu}} \\
\omega^i_{\barj} 
\end{pmatrix}
=
\begin{pmatrix}
t^{\nu}_{\bar{\mu} k} &
t^{\nu}_{\bar{\mu} \bar{k}} &
t^{\nu}_{\bar{\mu} \sigma} &
t^{\nu}_{\bar{\mu} \bar{\sigma}} &
t^{\nu \sigma}_{\bar{\mu} k} &
t^{\nu \bar{\sigma}}_{\bar{\mu} \bar{k}} \\
t^i_{\barj k} &
t^i_{\barj \bar{k}} &
t^i_{\barj \sigma} &
t^i_{\barj \bar{\sigma}} &
t^{i \sigma}_{\barj k} &
t^{i \bar{\sigma}}_{\barj \bar{k}}
\end{pmatrix}
\begin{pmatrix}
\omega^k \\
\omega^{\bar{k}} \\
\omega^{\sigma} \\
\omega^{\bar{\sigma}} \\
\omega^k_{\sigma} \\
\omega^{\bar{k}}_{\bar{\sigma}}
\end{pmatrix}
\]
so that
\begin{equation} \label{eqn:Middle}
\begin{pmatrix}
\omega^i_{\bar{\mu} j} \\
\omega^i_{\bar{\mu} \barj} \\
\omega^i_{\bar{\mu} \nu} \\
\omega^i_{\bar{\mu} \bar{\nu}} \\
\omega^{\nu}_{\bar{\mu}} \\
\omega^i_{\barj}
\end{pmatrix}
=
\begin{pmatrix}
t^i_{\bar{\mu} jk} &
t^i_{\bar{\mu} j \bar{k}} &
t^i_{\bar{\mu} j \sigma} &
t^i_{\bar{\mu} j \bar{\sigma}} &
- \delta^i_k t^{\sigma}_{\bar{\mu} j} &
t^i_{\bark j} \delta^{\bar{\sigma}}_{\bar{\mu}} \\
t^i_{\bar{\mu} \barj k} &
t^i_{\bar{\mu} \barj \bar{k}} &
t^i_{\bar{\mu} \barj \sigma} &
t^i_{\bar{\mu} \barj \bar{\sigma}} &
- \delta^i_k t^{\sigma}_{\bar{\mu} \barj} &
t^i_{\bar{k} \barj} \delta^{\bar{\sigma}}_{\bar{\mu}} \\
t^i_{\bar{\mu} \nu k} &
t^i_{\bar{\mu} \nu \bar{k}} &
t^i_{\bar{\mu} \nu \sigma} &
t^i_{\bar{\mu} \nu \bar{\sigma}} &
- \delta^i_k t^{\sigma}_{\barmu \nu} &
t^i_{\bark \nu} \delta^{\barsigma}_{\barmu} \\
t^i_{\bar{\mu} \bar{\nu} k} &
t^i_{\bar{\mu} \bar{\nu} \bar{k}} &
t^i_{\bar{\mu} \bar{\nu} \sigma} &
t^i_{\bar{\mu} \bar{\nu} \bar{\sigma}} &
- \delta^i_k t^{\sigma}_{\barmu \barnu} &
t^i_{\bark \barnu} \delta^{\barsigma}_{\barmu} \\
t^{\nu}_{\bar{\mu} k} &
t^{\nu}_{\bar{\mu} \bar{k}} &
t^{\nu}_{\bar{\mu} \sigma} &
t^{\nu}_{\bar{\mu} \bar{\sigma}} &
t^{\nu \sigma}_{\bar{\mu} k} &
t^{\nu \bar{\sigma}}_{\bar{\mu} \bar{k}} \\
t^i_{\barj k} &
t^i_{\barj \bar{k}} &
t^i_{\barj \sigma} &
t^i_{\barj \bar{\sigma}} &
t^{i \sigma}_{\barj k} &
t^{i \bar{\sigma}}_{\barj \bar{k}}
\end{pmatrix}
\begin{pmatrix}
\omega^k \\
\omega^{\bar{k}} \\
\omega^{\sigma} \\
\omega^{\bar{\sigma}} \\
\omega^k_{\sigma} \\
\omega^{\bar{k}}_{\bar{\sigma}}
\end{pmatrix}.
\end{equation}
Out of these equations, we can also see
that 
\begin{xalignat}{2}
t^{\sigma}_{\barmu \barnu} &= t^{\sigma}_{\barnu \barmu} &
t^i_{\bark \barnu} \delta^{\barsigma}_{\barmu}
&=  t^i_{\bark \barmu} \delta^{\barsigma}_{\barnu}.
\end{xalignat}
Plugging back into the equation \vref{eqn:POPO}
we find that
\begin{xalignat}{3} \label{eqn:Bloody}
\delta^i_j t^{\nu \sigma}_{\barmu k} &=
\delta^i_k t^{\sigma \nu}_{\barmu j} &
t^{i \barsigma}_{\barj \bark} \delta^{\barnu}_{\barmu}
&=
t^{i \barnu}_{\bark \barj} \delta^{\barsigma}_{\barmu} &
\delta^i_j t^{\nu \barsigma}_{\barmu \bark}
=
- t^{i \nu}_{\bark j} \delta^{\barsigma}_{\barmu}.
\end{xalignat}

Now we can take the exterior derivatives
of both sides of equation \vref{eqn:Middle}
and find out how the torsion transforms under 
$G_1$.}
Move up the fibers of $B_1$ to 
a subbundle $B_E$ on which
\begin{align*}
0 = t^{i \sigma}_{\barj i} = t^i_{\barj \sigma}.
\end{align*}
\extraStuff{
These equations, put together with equation
\vref{eqn:Bloody}, give us
\[
0 = t^i_{\barj \nu}
= t^{\nu \barsigma}_{\barmu \bark}
= t^{i \nu}_{\bark j}
\]
so that
\begin{equation} \label{eqn:Small}
\begin{pmatrix}
\omega^i_{\bar{\mu} j} \\
\omega^i_{\bar{\mu} \barj} \\
\omega^i_{\bar{\mu} \nu} \\
\omega^i_{\bar{\mu} \bar{\nu}} \\
\omega^{\nu}_{\bar{\mu}} \\
\omega^i_{\barj}
\end{pmatrix}
=
\begin{pmatrix}
t^i_{\bar{\mu} jk} &
t^i_{\bar{\mu} j \bar{k}} &
t^i_{\bar{\mu} j \sigma} &
t^i_{\bar{\mu} j \bar{\sigma}} &
- \delta^i_k t^{\sigma}_{\bar{\mu} j} &
t^i_{\bark j} \delta^{\bar{\sigma}}_{\bar{\mu}} \\
t^i_{\bar{\mu} \barj k} &
t^i_{\bar{\mu} \barj \bar{k}} &
t^i_{\bar{\mu} \barj \sigma} &
t^i_{\bar{\mu} \barj \bar{\sigma}} &
- \delta^i_k t^{\sigma}_{\bar{\mu} \barj} &
t^i_{\bar{k} \barj} \delta^{\bar{\sigma}}_{\bar{\mu}} \\
t^i_{\bar{\mu} \nu k} &
t^i_{\bar{\mu} \nu \bar{k}} &
t^i_{\bar{\mu} \nu \sigma} &
t^i_{\bar{\mu} \nu \bar{\sigma}} &
- \delta^i_k t^{\sigma}_{\barmu \nu} &
0 \\
t^i_{\bar{\mu} \bar{\nu} k} &
t^i_{\bar{\mu} \bar{\nu} \bar{k}} &
t^i_{\bar{\mu} \bar{\nu} \sigma} &
t^i_{\bar{\mu} \bar{\nu} \bar{\sigma}} &
- \delta^i_k t^{\sigma}_{\barmu \barnu} &
t^i_{\bark \barnu} \delta^{\barsigma}_{\barmu} \\
t^{\nu}_{\bar{\mu} k} &
t^{\nu}_{\bar{\mu} \bar{k}} &
t^{\nu}_{\bar{\mu} \sigma} &
t^{\nu}_{\bar{\mu} \bar{\sigma}} &
t^{\nu \sigma}_{\bar{\mu} k} &
0 \\
t^i_{\barj k} &
t^i_{\barj \bar{k}} &
0 &
t^i_{\barj \bar{\sigma}} &
0 &
t^{i \bar{\sigma}}_{\barj \bar{k}}
\end{pmatrix}
\begin{pmatrix}
\omega^k \\
\omega^{\bar{k}} \\
\omega^{\sigma} \\
\omega^{\bar{\sigma}} \\
\omega^k_{\sigma} \\
\omega^{\bar{k}}_{\bar{\sigma}}
\end{pmatrix}.
\end{equation}
}
This subbundle $B_E$ is right principal for
the subgroup $G_{\C{}} \subset G_1$. 
\extraStuff{
Indeed 
$G_{\C{}}$ consists
of the matrices of the form
\[
\begin{pmatrix}
a^i_j & 0 & 0 & 0 & 0 & 0 & 0 & 0 \\
0 & a^{\bari}_{\barj} & 0 & 0 & 0 & 0 & 0 & 0 \\
a^{\mu}_j & 0 & a^{\mu}_{\nu} & 0 
& 0 & 0 & 0 & 0 \\
0 & a^{\bar{\mu}}_{\barj} & 
0 & a^{\bar{\mu}}_{\bar{\nu}} & 
0 & 0 & 0 & 0 
\\
a^{i}_{\mu j} & 0 &
a^{i}_{\mu \nu} & 0 &
a^{i}_{j} A^{\nu}_{\mu} & 0 & 
0 & 0  
\\
0 & a^{\bari}_{\bar{\mu} \barj} &
0 & a^{\bari}_{\bar{\mu} \bar{\nu}} &
0 & 
a^{\bari}_{\barj} A^{\bar{\nu}}_{\bar{\mu}} &
0 & 
0
\\
0 & 0 &
0 & 0 &
0 & 
0 &
a^{i}_{j} A^{\bar{\nu}}_{\bar{\mu}} & 
0
\\
0 & 0 &
0 & 0 &
0 & 
0 &
0 & 
a^{\bari}_{\barj} A^{\nu}_{\mu}
\end{pmatrix}.
\]
}
This is just the same group $G_{\C{}}$ that occurs in
the $G_{\C{}}$ structure on the Grassmann
bundle of a  complex manifold.
It is clear that we can not reduce the structure
group any further in general, because in the
flat case (of a complex manifold) the biholomorphism
group acts transitively on this bundle $B_E = B_{\C{}}$.
\extraStuff{
Absorbing the torsion, we can arrange
\begin{equation} \label{eqn:TINY}
\begin{pmatrix}
\omega^i_{\barmu} \\
\omega^i_{\bar{\mu} j} \\
\omega^i_{\bar{\mu} \barj} \\
\omega^i_{\bar{\mu} \nu} \\
\omega^i_{\bar{\mu} \bar{\nu}} \\
\omega^{\nu}_{\bar{\mu}} \\
\omega^i_{\barj}
\end{pmatrix}
=
\begin{pmatrix}
0 &
0 &
0 &
0 &
0 &
0 \\
0 &
0 &
0 &
0 &
0 &
0  \\
0 &
0 &
0 &
0 &
0 &
t^i_{\bar{k} \barj} \delta^{\bar{\sigma}}_{\bar{\mu}} \\
0 &
0 &
0 &
0 &
0 &
0 \\
0 & 0 & 0 & 0 &
0 &
t^i_{\bark \barnu} \delta^{\barsigma}_{\barmu} \\
0 & 0 & 0 & 
0 &
t^{\nu \sigma}_{\bar{\mu} k} &
0 \\
0 &
t^i_{\barj \bar{k}} &
0 &
t^i_{\barj \bar{\sigma}} &
0 &
t^{i \bar{\sigma}}_{\barj \bar{k}}
\end{pmatrix}
\begin{pmatrix}
\omega^k \\
\omega^{\bar{k}} \\
\omega^{\sigma} \\
\omega^{\bar{\sigma}} \\
\omega^k_{\sigma} \\
\omega^{\bar{k}}_{\bar{\sigma}}
\end{pmatrix}
\end{equation}
with
\begin{xalignat*}{3}
t^i_{\barj \bark} &= -t^i_{\bark \barj} &
t^i_{\bark \barnu} \delta^{\barsigma}_{\barmu} &= 
t^i_{\bark \barmu} \delta^{\barsigma}_{\barnu} &
\delta^i_j t^{\nu \sigma}_{\barmu k} &=
\delta^i_k t^{\sigma \nu}_{\barmu j} \\
&& t^{i \barsigma}_{\barj \bark} \delta^{\barnu}_{\barmu} &=
t^{i \barnu}_{\bark \barj} \delta^{\barsigma}_{\barmu}. &&
\end{xalignat*}

Differentiating:\footnote{More precisely, 
taking exterior derivative modulo 
$\omega^j, \omega^{\barj}, \omega^{\barmu}, \omega^{\barj}_{\barmu}$
of the equation for $\omega^i_{\barj}$ in 
equation \vref{eqn:TINY}.}
\[
\begin{pmatrix}
\omega^{\mu}_{\barj} \\
\omega^i_{\mu \barj} 
\end{pmatrix}
=
\begin{pmatrix}
t^{\mu}_{\barj k} &
t^{\mu}_{\barj \bark} &
t^{\mu}_{\barj \sigma} &
t^{\mu}_{\barj \barsigma} &
t^{\mu \sigma}_{\barj k} &
t^{\mu \barsigma}_{\barj \bark} \\
t^{i}_{\mu \barj k} &
t^{i}_{\mu \barj \bark} &
t^{i}_{\mu \barj \sigma} &
t^{i}_{\mu \barj \barsigma} &
- \delta^i_k t^{\sigma}_{\barj \mu} &
t^{i \barsigma}_{\mu \barj \bark}
\end{pmatrix}
\begin{pmatrix}
\omega^k \\
\omega^{\bark} \\
\omega^{\sigma} \\
\omega^{\barsigma} \\
\omega^k_{\sigma} \\
\omega^{\bark}_{\barsigma}
\end{pmatrix}
\]
with
\[
\delta^i_m t^{\mu \sigma}_{\barj k} =
\delta^i_k t^{\sigma \mu}_{\barj m} \text{ and }
t^i_{\mu \barj \sigma} = t^i_{\sigma \barj \mu}.
\]

Absorbing torsion, we can
arrange 
\[
\begin{pmatrix}
\omega^{\mu}_{\barj} \\
\omega^i_{\mu \barj} 
\end{pmatrix}
=
\begin{pmatrix}
0 &
t^{\mu}_{\barj \bark} &
0 &
t^{\mu}_{\barj \barsigma} &
t^{\mu \sigma}_{\barj k} &
t^{\mu \barsigma}_{\barj \bark} \\
0 &
t^{i}_{\mu \barj \bark} &
0 &
t^{i}_{\mu \barj \barsigma} &
0 &
t^{i \barsigma}_{\mu \barj \bark}
\end{pmatrix}
\begin{pmatrix}
\omega^k \\
\omega^{\bark} \\
\omega^{\sigma} \\
\omega^{\barsigma} \\
\omega^k_{\sigma} \\
\omega^{\bark}_{\barsigma}
\end{pmatrix}
\]
and
\begin{xalignat*}{2}
t^{\mu}_{\barj \bark} &= - t^{\mu}_{\bark \barj} &
t^i_{\mu \barj \bark} &= - t^i_{\mu \bark \barj}.
\end{xalignat*}

Differentiating equation \vref{eqn:TINY}
we find\footnote{More precisely,
first we differentiate the equation for 
$\omega^{\nu}_{\barmu}$ in equation \vref{eqn:TINY},
and consider the result modulo $\omega^j,\omega^{\barj},
\omega^{\nu},\omega^{\barnu},\omega^{j}_{\nu}$.
That gives the first equation.
Then we differentiate the equation for $\omega^i_{\barmu \nu}$,
and consider the result modulo $\omega^j,\omega^{\barj},
\omega^{\sigma},\omega^{\barsigma},\omega^i_{\sigma}$,
to get the second and third equations.
Finally differentiate the equation for $\omega^i_{\barj}$
and consider the result modulo $\omega^k,\omega^{\bark},
\omega^{\sigma},\omega^{\barsigma},\omega^{\bark}_{\barsigma}$.
}
\begin{align*}
t^{\nu \barsigma}_{\barj \bark} \delta^{\bartau}_{\barmu}
&=
t^{\nu \bartau}_{\bark \barj} \delta^{\barsigma}_{\barmu} \\
t^{\nu \tau}_{\barmu k} = 0 &\text{ unless } d=1 \\
\left ( 
t^{i \bartau}_{\nu \barj \bark}
-
t^i_{\barj \barsigma} t^{\barsigma \bartau}_{\nu \bark}
\right ) \delta^{\barepsilon}_{\barmu}
&=
\left ( 
t^{i \barepsilon}_{\nu \bark \barj}
-
t^i_{\bark \barsigma} t^{\barsigma \barepsilon}_{\nu \barj}
\right ) \delta^{\bartau}_{\barmu} \\
\delta^i_k t^{\mu \nu}_{\barj m} 
&=
\delta^i_m t^{\nu \mu}_{\barj k}.
\end{align*}
(Recall that $d$ is the number of dependent
variables $w$ in the equation \vref{eqn:Eqn}.)
}

Absorbing torsion, our structure equations are:
\begin{equation}
\begin{split}
d 
\begin{pmatrix}
\omega^i \\
\omega^{\mu} \\
\omega^i_{\mu}
\end{pmatrix}
=&
-
\begin{pmatrix}
\omega^i_j & \omega^i_{\nu} & 0 \\
\omega^{\mu}_j & \omega^{\mu}_{\nu} & 0 \\
\omega^i_{\mu j} & \omega^i_{\mu \nu} &  \omega^i_j \delta^{\nu}_{\mu}
- \delta^i_j \omega^{\nu}_{\mu}
\end{pmatrix}
\wedge 
\begin{pmatrix}
\omega^j \\
\omega^{\nu} \\
\omega^j_{\nu}
\end{pmatrix}
\\
&-
\begin{pmatrix}
0 &  0 \\
t^{\mu \sigma}_{\barj k} & t^{\mu \sigma}_{\barnu k} \\
0 & 0 
\end{pmatrix}
\omega^k_{\sigma}
\wedge
\begin{pmatrix}
\omega^{\barj} \\
\omega^{\barnu}
\end{pmatrix}
\\&-
\begin{pmatrix}
t^i_{\barj \bark} &  
t^i_{\barj \barsigma} & 
t^{i \barsigma}_{\barj \bark}  \\
t^{\mu}_{\barj \bark} 
&  t^{\mu}_{\barj \barsigma} & 
t^{\mu \barsigma}_{\barj \bark} \\
t^i_{\mu \barj \bark} &  
t^i_{\mu \barj \barsigma} &  
t^{i \barsigma}_{\mu \barj \bark} 
\end{pmatrix}
\begin{pmatrix}
\omega^{\bark} \\
\omega^{\barsigma} \\
\omega^{\bark}_{\barsigma} 
\end{pmatrix}
\wedge \omega^{\barj}
\end{split}
\end{equation}
with
\begin{equation} \label{eqn:TorsionRelations} 
\begin{split}
t^i_{\barj \bark} &= -t^i_{\bark \barj} \\
t^{\mu}_{\barj \bark} &= - t^{\mu}_{\bark \barj} \\
t^i_{\mu \barj \bark} &= - t^i_{\mu \bark \barj} \\
t^i_{\bark \barnu} \delta^{\barsigma}_{\barmu} &= 
t^i_{\bark \barmu} \delta^{\barsigma}_{\barnu} \\
t^{\nu \barsigma}_{\barj \bark} \delta^{\bartau}_{\barmu}
&=
t^{\nu \bartau}_{\bark \barj} \delta^{\barsigma}_{\barmu} \\
t^{i \barsigma}_{\barj \bark} \delta^{\barnu}_{\barmu} &=
t^{i \barnu}_{\bark \barj} \delta^{\barsigma}_{\barmu} \\
\left ( 
t^{i \bartau}_{\nu \barj \bark}
-
t^i_{\barj \barsigma} t^{\barsigma \bartau}_{\nu \bark}
\right ) \delta^{\barepsilon}_{\barmu}
&=
\left ( 
t^{i \barepsilon}_{\nu \bark \barj}
-
t^i_{\bark \barsigma} t^{\barsigma \barepsilon}_{\nu \barj}
\right ) \delta^{\bartau}_{\barmu} \\
\delta^i_k t^{\mu \nu}_{\barj m} 
&=
\delta^i_m t^{\nu \mu}_{\barj k} \\
\delta^i_j t^{\nu \sigma}_{\barmu k} &=
\delta^i_k t^{\sigma \nu}_{\barmu j}. \\
\end{split}
\end{equation}

We have organized the structure equations into
$(2,0)+(1,1)$ and $(0,2)$ forms. Note that the
$\omega^i_j$ etc. in the first term might not
be $(1,0)$ forms, and therefore the first
term might contribute $(2,0)+(1,1)$ quantities.
Since the structure group acts in a complex representation,
it preserves an almost complex structure,
which is integrable (i.e. a complex structure)
exactly when the invariants in the $(0,2)$
part vanish. \extraStuff{This is just the equation
\[
0
=
\begin{pmatrix}
t^i_{\barj \bark} &  
t^i_{\barj \barsigma} & 
t^{i \barsigma}_{\barj \bark}  \\
t^{\mu}_{\barj \bark} 
&  t^{\mu}_{\barj \barsigma} & 
t^{\mu \barsigma}_{\barj \bark} \\
t^i_{\mu \barj \bark} &  
t^i_{\mu \barj \barsigma} &  
t^{i \barsigma}_{\mu \barj \bark} 
\end{pmatrix}.
\]
}

\subsection{Immediate corollaries of the structure equations}

\begin{proposition} Suppose that $E \subset \Gro{2n}{TM}$
is a system of differential
equations with Cauchy--Riemann tableau.
Then each point $P \in E$ in the fiber above a 
point $m \in M$ determines a
complex structure on $T_m M$; call it $J_P$. 
Moreover the plane $P \subset T_m M$ is $J_P$
complex linear.
\end{proposition}
\begin{proof}
The structure equations show that each coframing
in $B_E$ is determined up to complex linear
multiples. Moreover, the elements $\eta^i,\eta^{\mu}$
of an adapted coframe from $B_E$ are basic for the projection
to $M$, so form a coframe on $M$ at $m$. The
1-forms $\eta^i$ vanish on $P$. Therefore the
$\eta^i,\eta^{\mu}$ identify $T_m M$ with 
$\C{n+d}$, and identify $P$ with $\C{n}\oplus 0$,
a complex subspace. These $\eta^i$
are determined up to complex multiples,
as are the $\eta^i,\eta^{\mu}$ together.
Therefore all choices of coframes from
the fibers of $B_E$ determine the same almost
complex structure $J_P$ on $T_m M$.
\end{proof}

\begin{proposition}
Suppose that $E \subset \Gro{2n}{TM}$
is a system of differential
equations with Cauchy--Riemann tableau.
Then every $C^2$ solution of this system,
i.e. every immersion $f : \Sigma \to M$
whose tangent spaces belong to $E$,
is endowed with the structure of
a complex manifold.
\end{proposition}
\begin{proof}
Given any immersion $f : \Sigma \to M$
of an oriented manifold $\Sigma$ of dimension $2n$,
we can let $\hat{f} : \Sigma \to \Gro{2n}{TM}$
be the map associating to each point
of $\Sigma$ the tangent space
\[
\hat{f}(s) = f'(s) \cdot T_s \Sigma \in 
\Gro{2n}{TM}.
\]
To have $f$ solve $E$ means precisely
that $\hat{f}(s) \in E$ for all $s \in \Sigma$.
Clearly we get the diagram
\[
\xymatrix{
                                    & E \ar[d]^{\pi} \\
\Sigma \ar[r]^{f} \ar[ur]^{\hat{f}} & M 
}
\]
which differentiates to show that 
\[
\pi' \left ( \hat{f}(s) \right ) \hat{f}'(s) =
f'(s) : T_s \Sigma \to f'(s) T_s \Sigma = \hat{f}(s) \subset T_{f(s)} M
\]
so that
\[
\hat{f}'(s) T_s \Sigma \in \Theta \left ( \hat{f}(s) \right ).
\]
It also shows that \(  \hat{f}'(s) T_s \Sigma \) intersects
the tangent space to the fiber $E \to M$
transversely.

Take any local adapted coframing from
the bundle $B_E$, say $\eta^i,\eta^{\mu},\eta^i_{\mu}$,
and pull it back via $\hat{f}$.
We find that  $\eta^i=0$, and 
that $\eta^{\mu}$ is a coframing
on $\Sigma$ (actually, only on an open 
subset of $\Sigma$, since it is only a local
coframing), because of the transversality
of $\hat{f}$ with the fibers. 
But then on $\Sigma$ we have
\[
0 = d \eta^i = - \eta^i_{\mu} \wedge 
\eta^{\mu} 
\]
so that 
\[
\eta^i_{\mu} = P^i_{\mu \nu} \eta^{\nu}
\]
for some complex valued functions $P^i_{\mu \nu} = P^i_{\nu \mu}$ 
on $\Sigma$.
By the structure equations,
\[
d \eta^{\mu} = - \left ( \eta^{\mu}_{\nu}
- t^{\mu \sigma}_{\bartau k} P^k_{\sigma \nu} 
\eta^{\bartau} \right )
\wedge \eta^{\nu}.
\]
We see that there are no $(0,2)$ terms
appearing, and therefore (by the
Newlander--Nirenberg theorem) % , see
% Malgrange \cite{Malgrange:1972a,Malgrange:1972b}) 
these $\eta^{\mu}$
define a complex coframing giving a complex
structure on $\Sigma$. If we change the 
choice of coframing, then we obtain the 
same complex structure, because the structure
group acts via a complex representation.
\end{proof}

\begin{proposition} The manifold $E$ 
which parameterizes a differential
equation $E \to \Gro{2n}{TM}$
with Cauchy--Riemann tableau is
endowed invariantly with the structure
of an almost complex manifold, so that 
the stalks $E \to M$ are complex submanifolds,
and the planes $\Theta(P) = \pi'(P)^{-1}P$
are complex planes. 
\end{proposition}
\begin{proof} 
The tangent spaces to the 
stalks are described in terms  of any
adapted coframing $\eta$ by the 
equations
\[
\eta^i = \eta^{\mu} = 0
\] 
which are complex linear, so the
stalks are almost complex submanifolds.

The $\Theta$ planes are 
described (invariantly) by the equation
$\eta^i=0$ for any adapted coframing $\eta^i$
on $E$. These are complex linear equations,
so the $\Theta$ planes are complex planes.

To see that the stalks are actually complex
manifolds, plug in the structure equations
to see that
\[
d \eta^i_{\mu} = - \left ( \eta^i_j \delta^{\nu}_{\mu}
- \delta^i_j \eta^{\nu}_{\mu} \right ) \wedge \eta^j_{\nu}
\]
on these stalks. Therefore by
the Newlander--Nirenberg theorem,
the stalks are complex submanifolds.
\end{proof}

There are a host of bundles invariantly
defined on $E$, and we can read them off
of the structure equations. We will say
that a vector bundle $W \to E$ is \emph{soldered
by a representation} $\rho : G \to \GL{V}$
if $W$ is equipped with an isomorphism
with the vector bundle $(B \times V)/G$,
where $G$ acts via the diagonal action 
on $B \times V$. A list of some of
these vector bundles and their solderings
is given in table~\vref{tbl:solder}.
The complex vector bundle $\Xi$ is called the
\emph{characteristic vector bundle}.
In particular, it makes clear that
\begin{lemma}
Let $\Vertical$ be the vector bundle on $E$ of
vertical vectors for the map $\pi : E \to M$,
i.e. $\Vertical = \ker \pi'$. Then
\[
\Theta/\Vertical = \Xi \otimes_{\C{}} \Vertical^*.
\]
\end{lemma}

\begin{table} \label{tbl:solder}
\[
\begin{array}{cl}
\toprule
\textit{Vector bundle} & \textit{Soldering representation} \\
\midrule
TE & 
\begin{pmatrix}
a^i_j & 0 & 0 \\
a^{\mu}_j & a^{\mu}_{\nu} & 0 \\
a^i_{\mu j} & a^i_{\mu \nu} & a^i_j A^{\nu}_{\mu}
\end{pmatrix}
\\ \addlinespace
%\hline
\Lambda^{1,0} E &
{}^t \begin{pmatrix}
a^i_j & 0 & 0 \\
a^{\mu}_j & a^{\mu}_{\nu} & 0 \\
a^i_{\mu j} & a^i_{\mu \nu} & a^i_j A^{\nu}_{\mu}
\end{pmatrix}^{-1}
\\ \addlinespace
%\hline
\Theta &
\begin{pmatrix}
a^{\mu}_{\nu} & 0 \\
a^i_{\mu \nu} & a^i_j A^{\nu}_{\mu}
\end{pmatrix}
\\ \addlinespace
%\hline
\Vertical = \ker \pi' &
\begin{pmatrix}
a^i_j A^{\nu}_{\mu}
\end{pmatrix}
\\ \addlinespace
%\hline
\Xi=TE/\Theta &
\begin{pmatrix}
a^i_j
\end{pmatrix} 
\\ \addlinespace
%\hline
\Xi^* &
\begin{pmatrix}
A^j_i
\end{pmatrix}
\\ \addlinespace
%\hline
\Vertical^* &
\begin{pmatrix}
A^j_i a^{\mu}_{\nu}
\end{pmatrix}
\\
%\hline
\bottomrule
\end{array}
\]
\caption{How various vector bundles are soldered. The bundle $\Vertical$
is the bundle of vertical vectors for $E \to M$.}
\end{table}

\begin{proposition} Every system
of equations $E \subset \Gro{2n}{TM}$ 
with Cauchy--Riemann
tableau can be approximated to first
order by the Cauchy--Riemann equations,
i.e. near each point, we can find coordinates
$z^{\mu},w$ with our point at the origin
of these coordinates, so that equation~\vref{eqn:Eqn}
has $F$ and the first derivatives of $F$ 
vanishing at the origin. Conversely, every 
system of equations of the form of \ref{eqn:Eqn}
has Cauchy--Riemann tableau precisely
when such coordinates exist near any point.
\end{proposition}
\begin{proof}
The invariants preventing
a Cauchy--Riemann tableau are first order,
so equations
with vanishing first derivatives in the
functions $F^i_{\barmu}$ must have
Cauchy--Riemann tableau.

Let us prove the other direction.
First, we consider taking our equation
and carrying out a simple change of coordinates.
By translating the coordinates, we can arrange
that the point $(z,w)$ we are interested in is the
origin of coordinates. By rotation in these variables,
we can arrange that any chosen tangent plane
which satisfies the equation at first order is
taken to the plane $dw=0$. Now we can take any
function $f \left ( z,\bar{z},w,\bar{w} \right )$ and change coordinates
to $W=w+f$. As long as $df=0$ at
the origin, this change of coordinates preserves
our conditions, and we can easily see how
it changes our system of equations.
\begin{align*}
\pd{W}{\bar{z}} =& \pd{w}{\bar{z}} + \pd{f}{\bar{z}} \\
                =& F \left ( z,\bar{z},W-f,\bar{W}-\bar{f}, 
\pd{W}{z}-\pd{f}{z}, 
\pd{\bar{W}}{\bar{z}}-\pd{\bar{f}}{\bar{z}} \right ) 
+ \pd{f}{\bar{z}}.
\end{align*}
The quadratic $f$ has no impact on the lowest order
terms inside $F$, so that its effect is felt only
in the $\partial{f}/\partial{\bar{z}}$ term. We can
use this term to wipe out the $w$ and $\bar{w}$
linear terms in $F$, and to wipe out $z$ linear terms.
But as for $\bar{z}$ linear terms, we can only wipe
out those which are symmetric: the $\bar{z}$ terms in
\[
\pd{F^{\barmu}}{z^{\barnu}} + \pd{F^{\barnu}}{z^{\barmu}}.
\]

Now we will begin a more abstract approach. Suppose that
we have two equations of the form of~\ref{eqn:Eqn},
say $E_0$ and $E_1$. Each is equipped with a 
foliation $E_j \to M_j$ with complex structures
on the fibers. 

% Equipping the base manifolds $M_j$
% (locally) with flat Euclidean metrics, we find that
% the equivalence problem to match up 
% the bundles $E_j$ by maps which are holomorphic
% on the fibers and isometries on the base becomes
% an elliptic equivalence problem, and torsion-free,
% so by Malgrange's theorems in Malgrange 
% \cite{Malgrange:1972a,Malgrange:1972b},
By the Newlander--Nirenberg theorem ``with
parameters'' (i.e. making obvious modifications
to the Malgrange proof which is clearly presented
in Nirenberg \cite{Nirenberg:1973}, to allow
for families of complex structures, using
smooth dependence of solutions of elliptic
equations on parameters),
locally we can find a diffeomorphism
\[
\xymatrix{
E_0 \ar[r] \ar[d] & E_1 \ar[d] \\
M_0 \ar[r] & M_1
}
\]
which is holomorphic on the fibers of the $E_j$.
It is easy to see (in local coordinates)
that the pseudogroup of diffeomorphisms
of $E_0$ which are holomorphic on the fibers
acts transitively on coframes of $E_0$
of the form $\eta^i,\eta^{\mu},\eta^i_{\mu}$
with the $\eta^i_{\mu}$ being $(1,0)$-forms
on the fibers, and the $\eta^i,\eta^{\mu}$
vanishing on the fibers. Therefore we
can arrange that some local adapted 
coframes $\eta$ of $E_0$  and $\xi$ of $E_1$ agree
at some point $P \in E_0$. By perhaps altering
the choices of the coframes, they
must therefore satisfy equations of the
form
\[
\begin{pmatrix}
\eta^i \\
\eta^{\bari} \\
\eta^{\mu} \\
\eta^{\barmu} \\
\eta^i_{\mu} \\
\eta^{\bari}_{\barmu}
\end{pmatrix}
=
\begin{pmatrix}
\delta^i_j & a^i_{\barj} & 0 & a^i_{\barnu} & 0 & 0 \\
a^{\bari}_j & \delta^{\bari}_{\barj} & a^{\bari}_{\nu} & 0 & 0 & 0 \\
0 & a^{\mu}_{\barj} & \delta^{\mu}_{\nu} & a^{\mu}_{\barnu} & 0 & 0 \\
a^{\barmu}_j & 0 & a^{\barmu}_{\nu} & \delta^{\barmu}_{\barnu} & 0 & 0 \\ 
0 & a^i_{\mu \barj} & 0 & a^i_{\mu \barnu} & \delta^i_j \delta^{\nu}_{\mu} + a^{i \nu}_{j \mu} & 0 \\
a^{\bari}_{\barmu j} & 0 & a^{\bari}_{\barmu \nu} & 0 & 0 & \delta^{\bari}_{\barj} \delta^{\barnu}_{\barmu}
+ a^{\bari \barnu}_{\barj \barmu} 
\end{pmatrix}
\begin{pmatrix}
\xi^j \\
\xi^{\barj} \\
\xi^{\nu} \\
\xi^{\barnu} \\
\xi^j_{\nu} \\
\xi^{\barj}_{\barnu}
\end{pmatrix}.
\]
(The last two rows follow from differentiating
the first four.) All of the functions $a$ 
vanish at the chosen point $P$. Differentiating, we find
\[
\begin{pmatrix}
\eta^i_j - \xi^i_j \\
da^i_{\barj} \\
da^i_{\barnu}
\end{pmatrix}
=
\begin{pmatrix}
a^i_{jk} & a^i_{j \bark} & a^i_{j \barsigma} \\
a^i_{\barj k} & a^i_{\barj \bark} & a^i_{\barj \barsigma} \\
a^i_{\barnu k} & a^i_{\barnu \bark} & a^i_{\barnu \barsigma} 
\end{pmatrix}
\begin{pmatrix}
\xi^k \\
\xi^{\bark} \\
\xi^{\barsigma}
\end{pmatrix}
\]
with symmetries in the lower indices from Cartan's lemma.
If we take the equation $E_1$ to be a complex
structure, i.e. the Cauchy--Riemann equations,
then we can easily see from these structure
equations that near the point $P$,
the numbers $a^i_{\barnu k}, a^i_{\barnu \bark}, a^i_{\barnu \barsigma}$
correspond to the numbers
\[
\pd{F^i_{\barnu}}{w^k}, 
\pd{F^i_{\barnu}}{w^{\bark}}, 
\pd{F^i_{\barnu}}{z^{\barmu}}.
\]
But by explicit coordinate manipulations above, we managed
to kill these terms.
\end{proof}
\section{High dimension and codimension}
Suppose that $n > 1$ and $d > 1$, i.e. 
the equation \vref{eqn:Eqn} has
more than one
(independent) $z$ variable and more than
one (dependent) $w$ variable.
\extraStuff{
Using only the hypothesis that $n>1$, 
and taking the exterior derivative 
of the expression for $\omega^i_{\barj}$, 
modulo $\omega^k, \omega^{\bark}$,
we find
\[
0 = t^{\mu}_{\barj \barnu} =
t^i_{\mu \barj \barnu} 
= t^i_{\barj \bark}.
\]
Differentiating the equation $\omega^{\mu}_{\barnu}=0$,
we obtain the equation
\[
0 = \omega^{\mu}_{\barnu \barj} 
- t^{\mu}_{\bark \barj} \omega^{\bark}_{\barnu} 
\pmod{\omega^{m},\omega^{\barm},\omega^{\tau},
\omega^{\bartau}}.
\]
Differentiating the equation 
$\omega^{\mu}_{\barj}=t^{\mu}_{\barj \bark} \omega^{\bark}
+ t^{\mu \barsigma}_{\barj \bark} \omega^{\bark}_{\barsigma}$,
we obtain the equation
\[
0 = \omega^{\mu}_{\barnu \barj} 
+ t^{\mu}_{\bark \barj} \omega^{\bark}_{\barnu} 
\pmod{\omega^{m},\omega^{\barm},\omega^{\tau},
\omega^{\bartau}}.
\]
Putting these two together, we find
\[
t^{\mu}_{\bark \barj} = 0.
\]
Finally, differentiating the equation for $d \omega^i$,
and writing the result modulo $\omega^{j},\omega^{\sigma} \wedge 
\omega^{\nu},\omega^{j}_{\nu}$,
we find that 
\[
t^i_{\mu \barj \bark} = 0.
\]
}
The  torsion equations 
force all of the torsion coefficients
to vanish. By the Cartan--K{\"a}hler theorem
every real analytic system of partial
differential equations with Cauchy--Riemann
tableau for $d>1$ complex functions of 
$n>1$ complex variables becomes the Cauchy--Riemann
equations in some system of local coordinates.
In fact, real analyticity is not needed:
\begin{theorem} \label{thm:FFFF}
Every system of partial 
differential equations of the form 
\vref{eqn:Eqn} for $d>1$ complex functions
$w$ of $n>1$ complex variables $z$ 
which has Cauchy--Riemann tableau
becomes the Cauchy--Riemann equations
in some system of local coordinates.
\end{theorem}
\begin{proof} 
% This is an immediate
% consequence of Malgrange's proof of the flatness
% of involutive elliptic $G$ structures;
% see Malgrange \cite{Malgrange:1972a,Malgrange:1972b}.
% Alternatively, it is not difficult to prove this theorem 
Use the Newlander--Nirenberg theorem to produce
holomorphic coordinates on $E$.
\extraStuff{Follow the proof of proposition~\vref{prop:StrucEqns}.}
\end{proof}
\section{Small dimension and codimension}
A thorough study of these equations is
given in McKay \cite{McKay:1999,McKay:AlsoUnpub}.
\section{Hypersurfaces}
Suppose now that the number $n$ of independent
variables in the equation \vref{eqn:Eqn} is
greater than one, but that the number of 
dependent variables $d=1$, and that
the equation has Cauchy--Riemann tableau.  So this equation 
represents a 
generalization of the theory
of complex hypersurfaces in a complex manifold.
We will call it a \emph{hypersurface equation}.
\extraStuff{
Using the same computations as in the
previous section, we obtain the equations
\[
0 = t^{\mu}_{\barj \barnu} =
t^i_{\mu \barj \barnu} 
= t^i_{\barj \bark}.
\]
Differentiating the structure equations
we obtain the equations
\[
0 = t^{\mu}_{\barj \bark} = t^i_{\mu \barj \bark}
\]
and
\[
t^{\mu \nu}_{\barj k} = t^{\nu \mu}_{\barj k}
\quad \text{and} \quad
t^{\mu \nu}_{\bartau k} = t^{\nu \mu}_{\bartau k}.
\]
}
The structure equations are:
\begin{equation} \label{eqn:HypStruc}
d 
\begin{pmatrix}
\omega^i \\
\omega^{\mu} \\
\omega^i_{\mu}
\end{pmatrix}
=
-
\begin{pmatrix}
\omega^i_j & \omega^i_{\nu} & 0 \\
\omega^{\mu}_j & \omega^{\mu}_{\nu} & 0 \\
\omega^i_{\mu j} & \omega^i_{\mu \nu} &  \omega^i_j \delta^{\nu}_{\mu}
- \delta^i_j \omega^{\nu}_{\mu}
\end{pmatrix}
\wedge 
\begin{pmatrix}
\omega^j \\
\omega^{\nu} \\
\omega^j_{\nu}
\end{pmatrix}
-
\begin{pmatrix}
0 &  0 \\
t^{\mu \sigma}_{\barj k} & t^{\mu \sigma}_{\barnu k} \\
0 & 0 
\end{pmatrix}
\omega^k_{\sigma}
\wedge
\begin{pmatrix}
\omega^{\barj} \\
\omega^{\barnu}
\end{pmatrix}.
\end{equation}
The structure group preserves a complex
structure on $E \subset \Gro{2n}{TM}$, since
there are no $(0,2)$ terms in the structure
equations. 
\begin{proposition} For any choice of 
complex constants $T^{\mu \sigma}_{\barj k}, 
T^{\mu \sigma}_{\barnu k}$
satisfying
\[
T^{\mu \nu}_{\barj k} = T^{\nu \mu}_{\barj k}
\quad \text{and} \quad
T^{\mu \nu}_{\bartau k} = T^{\nu \mu}_{\bartau k}
\]
there is a 
hypersurface equation $E$, i.e. one of the form \vref{eqn:Eqn} with
Cauchy--Riemann tableau, for one complex
function $w$ of several complex variables $z^{\mu}$,
so that the associated $G$ structure $B_E$
satisfies the structure equations
\vref{eqn:HypStruc} so that
\[
T^{\mu \nu}_{\barj k} = t^{\mu \nu}_{\barj k}
\quad \text{and} \quad
T^{\mu \nu}_{\bartau k} = t^{\mu \nu}_{\bartau k}
\]
at some point of $B_E$. 
The general real analytic hypersurface
equation depends on $2$ real functions
of $2n+1$ real variables.
\end{proposition}
\begin{proof} This is immediate from the
Cartan--K{\"a}hler theorem, since the 
structure equations \vref{eqn:HypStruc}
are involutive.
\end{proof}

\subsection{Contact geometry}

So we see that such equations exist, 
and we wonder how to construct them.

\begin{proposition} The manifold $E \subset \Gro{2n}{TM}$
is a complex contact manifold, with contact plane
field $\Theta$, and the fibers of $\pi : E \to M$
are holomorphic Legendre submanifolds.
\end{proposition}
\begin{proof} Fattening up the structure
group to the complex contact group (see McKay \cite{McKay:Unpub})
we obtain the structure equations
\extraStuff{
\[
d 
\begin{pmatrix}
\omega^i \\
\omega^{\mu} \\
\omega^i_{\mu}
\end{pmatrix}
=
-
\begin{pmatrix}
\omega^i_j & \omega^i_{\nu} & 0 \\
\omega^{\mu}_j & \omega^{\mu}_{\nu} & \omega^{\mu \nu}_j \\
\omega^i_{\mu j} & \omega^i_{\mu \nu} &  \omega^i_j \delta^{\nu}_{\mu}
- \delta^i_j \omega^{\nu}_{\mu}
\end{pmatrix}
\wedge 
\begin{pmatrix}
\omega^j \\
\omega^{\nu} \\
\omega^j_{\nu}
\end{pmatrix}
\]
}
of a complex contact structure.
By the Newlander--Nirenberg theorem and the
Darboux theorem for complex contact structures,
% Malgrange \cite{Malgrange:1972a,Malgrange:1972b},
such structures are all locally isomorphic.
$\Theta$ is the contact
plane field. The fibers of $\pi : E \to M$
are given by the complex linear equations 
$\eta^i=\eta^{\mu}=0$, for any section $\eta$
of $B_E$. Therefore they are complex
submanifolds, and clearly tangent to the
plane field $\Theta$, therefore Legendre.
\end{proof}

Now we know how to construct equations
with Cauchy--Riemann tableau for one
complex function of several complex
variables: we simply construct a 
Legendre fibration, not necessarily
holomorphic, of a complex contact manifold $E$
of complex dimension $2n+1$.
Then the base of the fibration is 
our manifold $M$, and the 
manifold $E$ has a well defined 
immersion into $\Gro{2n}{TM}$.
It is in this sense that I say
that such equations exist microlocally:
locally on $E$.

Moreover, we can see that the equations
are locally solvable: solutions are just
Legendre submanifolds transverse to 
the fibers of the Legendre fibration.
Therefore local solutions exist, and depend
on $2$ real functions of $n$ real variables.
(Global solvability is more difficult;
we see this already in the context of complex manifolds
when trying to globally construct hypersurfaces.
For example, generic non-K\"ahler tori
have no compact hypersurfaces.)

Now we wish to consider the global
geometry of this immersion, or
in other words, to describe a
hypersurface equation locally on $M$.
\begin{theorem}[Merkulov \cite{Merkulov:1997}]
Let $E$ be a complex contact manifold
with contact plane field $\Theta$.
The (infinitely many) obstructions
to deforming a holomorphic Legendre
manifold $L$ of $E$
are found in the first cohomology
group $\Cohom{1}{\Xi_L}$ where $\Xi_L$ is the 
pullback to $L$ of the characteristic
line bundle $\Xi=TE / \Theta$. 
If these obstructions vanish then the
submanifold $L$ admits a locally
complete moduli space $\mathcal{M}$
of deformations
whose tangent space at $L$ is
$\Cohom{0}{\Xi_L}$. 
\end{theorem}
\extraStuff{
We will also use some elementary theory
of vector bundles on projective space.
\begin{lemma}
The cohomology groups of the line
bundles on $\CP{n}$ are
\[
\Cohom{p}{\mathcal{O}(k)}
=
\begin{cases}
0 & k < 0, \ p \ne n \\
0 & k \ge 0, \ p \ne 0 \\
\binom{n-k}{-k} & k \le 0, \ p=n \\
\binom{n+k}{k} & k \ge 0, \ p=0 
\end{cases}
\]
\end{lemma}
For proof see Griffiths \& Harris \cite{GriffithsHarris:1978}, pg. 156.

\begin{lemma} 
On $\CP{n}$
\[
\Cohom{0}{\mathcal{O}(1) \otimes_{\C{}} \Lambda^{1,0}} = 0
\]
where $\Lambda^{1,0}$ means the holomorphic
cotangent bundle of $\CP{n}$.
\end{lemma}
\begin{proof} Suppose that we have a holomorphic
section $s$ which does not vanish at some
point $p$. Take affine
coordinates near $p.$
Then in those coordinates
\[
s = \sigma \otimes a_j \, dz^j.
\]
By taking a linear change of those coordinates,
we can arrange
\[
s = s_1 = \sigma_1 \otimes dz^1.
\] 
Now by rotating those coordinates, we get
other sections
\[
s_j = \sigma_j \otimes dz^j
\]
not vanishing at $p.$
Wedging them together, we get a
section
\begin{align*}
\sigma_1 \dots \sigma_n \otimes dz^1 \wedge \dots \wedge dz^n
\in \operatorname{Det} \left (  
                \mathcal{O}(1) \otimes_{\C{}} \Lambda^{1,0} 
        \right )
&=
\mathcal{O}(1)^{\otimes n} \otimes_{\C{}}
\Lambda^{n,0}
\\
&=
\mathcal{O}(n) \otimes_{\C{}} \mathcal{O}(-(n+1))
\\
&=
\mathcal{O}(-1)
\end{align*}
not vanishing at $p.$
But this line bundle has no global sections,
because it has negative first Chern class.
\end{proof}
}
\begin{definition} By the term \emph{local geometry}
applied to a system of differential equations
$E \subset \Gro{2n}{TM}$ with Cauchy--Riemann
tableau, we mean the
three numbers
\begin{align*}
& \dim \Cohom{0}{\Xi_{E_m}} \\
& \dim \Cohom{1}{\Xi_{E_m}} \\
& \dim \Cohom{0}{\Xi_{E_m} \otimes \Lambda^{1,0}}
\end{align*}
where $\Xi$ is the quotient line bundle
$\Xi=TE / \Theta$, and $\Xi_{E_m}$ is the
pullback of that line bundle to a fiber
$E_m \subset E$, and $\Lambda^{1,0}$ is
the holomorphic cotangent bundle of $E_m$. 
These numbers are actually
integer valued functions on $M$.
\end{definition}
\extraStuff{
The rigidity of the complex
structure of complex projective space is 
proven in Fr{\"o}licher \& Nijenhuis 
\cite{FrolicherNijenhuis:1957}. 
An old conjecture, still unproven,
states that there is a unique complex
structure on complex projective space. We
do not need such a result.
}
\begin{theorem} \label{thm:BIG}
Let $E \subset \Gro{2n}{TM}$
be a system of equations of the form \vref{eqn:Eqn}
with Cauchy--Riemann tableau.
Suppose that $E$ has the same local
geometry as the Cauchy--Riemann equations.
Then $M$ bears a unique complex structure so that
$E$ is the set of complex hyperplanes in the
tangent spaces of $M$. In particular,
the system of equations \ref{eqn:Eqn} is
the Cauchy--Riemann equations in appropriate
coordinates.\footnote{It is in this sense that
I say that new theories of 
pseudoholomorphic hypersurfaces
do not exist locally, although they exist
microlocally.
By comparison, 
theorem~\vref{thm:FFFF} says that
new theories of high dimensional
and codimensional pseudoholomorphic
objects do not exist even microlocally.}
\end{theorem}
\begin{remark}
The reader will note that
if $E$ is the set of complex hyperplanes
in the tangent spaces of a complex
manifold $M$, i.e. $E$ is the Cauchy--Riemann
equations, then the fibers of $E \to M$
are copies of $\CP{n}$ and the line bundle
$\Xi_{E_m}$ is $\mathcal{O}(1)$.
For any $E_m$ fiber, the vector bundle $\Theta/TE_m$
is 
\[
\Theta/TE_m = \Xi \otimes_{\C{}} \Lambda^{1,0}
\]
where $\Lambda^{1,0}$ is the holomorphic 
cotangent bundle. So for the Cauchy--Riemann
equations, 
\begin{align*}
\dim_{\C{}} \Cohom{0}{\Xi_{E_m}} &= n+1 \\
\dim_{\C{}} \Cohom{1}{\Xi_{E_m}} &= 0 \\
\dim_{\C{}} \Cohom{0}{\Xi_{E_m} \otimes \Lambda^{1,0}}
&= 0
\end{align*}

Moreover local topology is unaffected
by small perturbations, i.e. under
deformation of the Cauchy--Riemann equations 
$E_0$ to equations $E_t$, the
rigidity of $\CP{n}$ and of its line 
bundles ensures that
all of the equations $E_t$ for $t$ near $0$
have the same local topology. This
holds true even if the equations $E_t$
are only immersed submanifolds of $\Gro{2n}{TM}$.
\end{remark}
\begin{proof} 
We have seen that the points of $M$
can be interpreted as Legendre submanifolds
in $E$: the fibers $E_m$. 
The cohomology numbers above
prove the existence of a locally
complete moduli space $\mathcal{M}$ 
of Legendre submanifolds, by Merkulov's
theorem.
Using a local section of $E \to M$,
we get a map of $M$ into the moduli
space $\mathcal{M}$. The spaces $M$ and $\mathcal{M}$
have the same
dimension, and $M$ is mapped by a smooth
injection. We want to show that this 
map is an immersion, and therefore
a local diffeomorphism. 

Consider a particular fiber $E_m$ of 
$E \to M$ over a point $m \in M$.
Let $\eta^i,\eta^{\mu},\eta^i_{\mu}$ 
be any section of $B_E$, i.e. any
adapted coframing. We will write it as
$\eta,\eta^{\mu},\eta_{\mu}$ since there is
only one value for $i$. As Merkulov \cite{Merkulov:1997}
explains, a tangent vector $v \in M$ 
corresponds to a unique section of
$\Cohom{0}{\Xi_{E_m}}$ which is determined
as follows: take the holomorphic section
of the normal bundle
\[
s^0 \pd{}{\eta} + s^{\mu} \pd{}{\eta^{\mu}}
\in \Cohom{0}{N_{E_m}}
\]
which projects to $v$. Then project it
to
\[
s^0 \pd{}{\eta} \in \Cohom{0}{\Xi_{E_m}}.
\]
If this vanishes, then 
the section of the
normal bundle must be the image of 
a section of $\Theta_{E_m}/TE_m$.
As we have seen
\[
\Theta_{E_m}(P) / T_P E_m
= P = \Xi(P) \otimes_{\C{}} \Lambda^{1,0}.
\] 

Again, this bundle has no global holomorphic
sections. 
Consequently, the vectors $v \in T_m M$
are injectively mapped to elements
of $\Cohom{0}{\Xi_{E_m}}=T_{E_m} \mathcal{M}$,
so the map $M \to \mathcal{M}$ is
a local diffeomorphism.

We can put a complex structure
on $M$, pulling back the one from $\mathcal{M}$,
and the map $E \to M$ must be holomorphic
for that complex structure, since the
map to moduli space $\mathcal{M}$ is. Choosing any
local holomorphic coordinates $w,z^{\mu}$
on $M$, we can pull them back to $E$ to
find that in terms of any adapted coframing
$\eta,\eta^{\mu},\eta_{\mu}$,
\begin{align*}
dw &= a \eta + a_{\nu} \eta^{\nu} \\
dz^{\mu} &= a^{\mu} \eta + a^{\mu}_{\nu} \eta^{\nu}
\end{align*}
for some complex valued functions $a,a_{\nu},a^{\mu},a^{\mu}_{\nu}$.
By a complex linear change of coordinates, we
can arrange that at some chosen point of $E$,
\[
\begin{pmatrix}
a & a_{\nu} \\
a^{\mu} & a^{\mu}_{\nu} 
\end{pmatrix}
= I.
\]
Taking exterior derivative,
we find
\[
t^{\mu \sigma}_{\barj k} = t^{\mu \sigma}_{\barnu k} = 0
\]
so that the $G$ structure is torsion-free,
and again by the Newlander--Nirenberg theorem,
and Darboux's theorem for complex contact
structures, 
% by Malgrange's theorem on torsion-free
% involutive elliptic structures, 
% Malgrange \cite{Malgrange:1972a,Malgrange:1972b},
it is flat. Therefore the equations \vref{eqn:Eqn}
are the Cauchy--Riemann equations.
\end{proof}

\begin{corollary} Every continuously varying 
family $E_t \subset \Gro{2n}{TM}$
of equations with Cauchy--Riemann tableau
for one complex function $w$ of several complex
variables $z$ with $E_0$ being the Cauchy--Riemann
equations
is a deformation of complex structures.
In other words, all $E_t$ are Cauchy--Riemann
equations.
\end{corollary}

\subsection{Geometry of the characteristic line bundle}

Let $E_m$ be a single fiber of $E \to M$.
Let $B_{E_m}$ be the quotient 
of the pullback
bundle of $B_E$ to $E_m \subset E$
by the group of matrices
\[
\begin{pmatrix}
1 & 0 & 0 & 0  \\
0 & 1 & 0 & 0  \\
a^i_{\mu j} & a^i_{\mu \nu} & 0 & 0 
\end{pmatrix}
\]
(which is a subgroup of the structure group of $B_E$).
Then on $B_{E_m}$ we have $\omega^i=\omega^{\mu}=0$.
But the 1-forms $\omega^i_{\mu}$ are 
semibasic for the projection $B_{E_m} \to E_m$.
The 1-forms $\omega^i_j,\omega^{\mu}_{\nu},
\omega^{\mu}_j$ are not uniquely defined
on $B_E$, being defined only up to adding
multiples of $\omega^i$ and $\omega^{\mu}$.
But on $B_{E_m}$ they are thereby uniquely
defined, since those vanish.  The 1-forms
$\omega^i_{\mu j}$ and $\omega^i_{\mu \nu}$ are
not well defined on $B_{E_m}$. The structure equations on $B_{E_m}$
are
\[
d
\begin{pmatrix}
\omega^i_j & \omega^i_{\nu} \\
\omega^{\mu}_j & \omega^{\mu}_{\nu}
\end{pmatrix}
=
-
\begin{pmatrix}
\omega^i_k & \omega^i_{\sigma} \\
\omega^{\mu}_k & \omega^{\mu}_{\sigma}
\end{pmatrix}
\wedge
\begin{pmatrix}
\omega^k_j & \omega^k_{\nu} \\
\omega^{\sigma}_j & \omega^{\sigma}_{\nu}
\end{pmatrix}
-
\begin{pmatrix}
0 & 0 \\
t^{\mu \tau}_{\barsigma m}
t^{\barsigma \barepsilon}_{j \barp}  &
t^{\mu \tau}_{\barsigma m}
t^{\barsigma \barepsilon}_{\nu \barp} 
\end{pmatrix}
\omega^m_{\tau} \wedge
\omega^{\barp}_{\barepsilon}.
\]
The bundle $B_{E_m} \to E_m$ is a principal right
$H$ bundle, where $H$ is the group
of complex matrices of the 
form 
\[
\begin{pmatrix}
a & 0 \\
b & c
\end{pmatrix}
\]
with $a$ a $1 \times 1$, $b$ an $n \times 1$ and $c$ an $n \times n$.
The elements of $B_{E_m}$ are identified with
the coframes $\eta^i,\eta^{\mu}$ which
belong to adapted coframes $\eta^i,\eta^{\mu},\eta^i_{\mu}$
from $B_E$ above the point $m \in M$.
Since $\eta^i,\eta^{\mu}$ are semibasic
for the projection $E \to M$, they
can be identified with a coframe on $M$
itself. This identifies $B_{E_m}$ with
a principal $H$ subbundle of the
$\GL{2n+2,\R{}}$ bundle of linear isomorphisms of $T_m M$
with $\R{2n+2}$ (the base of this bundle is a point).

Let us first consider the vector bundle over
$B_{E_m}$ whose fiber above a point 
$\left ( \eta^i, \eta^{\mu} \right ) \in B_{E_m}$
is just the vector space $T_m M$.
We will just call this bundle $E_m \times T_m M$. 
It is topologically trivial, since all of the
fibers are the same, but it has a complex structure as
well, given by using $\eta^i,\eta^{\mu}$ 
to identify $T_m M$ with $\C{n+1}$. 
This vector bundle is soldered by the
1-form
\[
\begin{pmatrix}
\omega^i_j & 0 \\
\omega^{\mu}_j & \omega^{\mu}_{\nu}
\end{pmatrix}.
\]
This gives it the structure of a complex
vector bundle.
The choice of a single constant vector $v \in T_m M$
gives rise to a section $\sigma_v$ of the
bundle $E_m \times T_m M$, 
represented as functions on $B_{E_m}$
given by
\[
F_v \left ( \eta^i,\eta^{\mu} \right )
=
\begin{pmatrix}
\eta^i(v) \\
\eta^{\bari}(v) \\
\eta^{\mu}(v) \\
\eta^{\barmu}(v)
\end{pmatrix}
\]
satisfying
\begin{equation}\label{eqn:DiffSigma}
d
\begin{pmatrix}
F^i \\
F^{\bari} \\
F^{\mu} \\
F^{\barmu}
\end{pmatrix}
=
-
\begin{pmatrix}
\omega^i_j & \omega^i_{\barj} & \omega^i_{\nu} & \omega^i_{\barnu} \\
\omega^{\bari}_j & \omega^{\bari}_{\barj} & \omega^{\bari}_{\nu} & \omega^{\bari}_{\barnu} \\
\omega^{\mu}_j & \omega^{\mu}_{\barj} & \omega^{\mu}_{\nu} & \omega^{\mu}_{\barnu} \\
\omega^{\barmu}_j & \omega^{\barmu}_{\barj} & \omega^{\barmu}_{\nu} & \omega^{\barmu}_{\barnu}
\end{pmatrix}
\begin{pmatrix}
F^j \\
F^{\barj} \\
F^{\nu} \\
F^{\barnu}
\end{pmatrix}.
\end{equation}
But in our situation, we see that on $B_{E_m}$
the 1-forms $\omega^i_{\barj}$ and $\omega^i_{\barnu}$ vanish.
Clearly these $F_v$ are not holomorphic sections of this
vector bundle, but if we quotient out the $F^{\mu}$
parts, the $F^i$ parts are holomorphic, sections
of the bundle $\Xi$ (whose fiber at $P$ is $T_m M/P$).
Indeed these are the holomorphic sections
of $\Xi$ we saw in the contact geometry above.

If we pick a point $P_0 \in E_m$, which
is a plane of codimension 2 in $T_m M$,
we can pick a $J_{P_0}$ complex basis $v_0, v_1, \dots, v_n$
for $T_m M$, with $v_1, \dots, v_n \in P_0$,
and try to construct a map to $\CP{n}$
by
\[
P \in E_m \mapsto 
\left [
\sigma_{v_0}(P) : \dots : \sigma_{v_n}(P)
\right ]
\in \CP{n}.
\]
This map is defined near $P_0$
because $\sigma_{v_0} \left ( P_0 \right ) \ne 0$
by construction.
It is easy to check from equation~\vref{eqn:DiffSigma}
that this map is an immersion near $P_0$.
One does this by taking $\eta^i,\eta^{\mu}$
which take $v_0, \dots, v_n$ to the standard
basis of $\C{n}$, and checking the differentials.

This would appear to determine a biholomorphism
$E_m \to \CP{n}$, but we must be very careful.
The basis $v_0, \dots, v_n$ is a complex
$J_{P_0}$ basis, but there might be a point
$P \in E_m$ where there are $J_P$ complex
linear relations among $v_0, \dots, v_n$.
For noncompact $E_m$ fibers this can occur.
At such points this map is not defined.
In fact, this map is clearly meromorphic, since
the set of such points is the set of zeros
of 
\[
\sigma_{v_0} \wedge \dots \wedge \sigma_{v_n}.
\]

\begin{lemma}
The map
\[
E_m \to \CP{}\left(\Cohom{0}{\Xi}^*\right)
\]
is well defined and a holomorphic immersion.
In particular, if $E_m$ is compact,
then $E_m$ is a smooth projective variety.
\end{lemma}
\begin{proof}
That this map is well defined follows
from there being, at each point $P \in E_m$,
some vector $v_0 \in T_m M \backslash P$.
That this map is an immersion follows
immediately from our discussion of the
map
\[
P \mapsto
\left [
\sigma_{v_0}(P) : \dots : \sigma_{v_n}(P)
\right ].
\]
\end{proof}

The central problem we face is that there
might be sections $\sigma \in \Cohom{0}{\Xi}$
which are not of the form $\sigma=\sigma_v$
for some vector $v \in T_m M$.
We will find some topological conditions
under which we can ensure that every
section of the characteristic
line bundle $\Xi$ has this
form.

\begin{lemma} \label{lm:One}
Suppose that every holomorphic section
$\sigma \in \Cohom{0}{\Xi}$ has the
form $\sigma=\sigma_v$ for some vector
$v \in V$. Suppose further that $E_m$ is
compact. Then $E_m$ is biholomorphic
to $\CP{n}$ via a biholomorphism
which identifies $\Xi$ with $\mathcal{O}(1)$.
\end{lemma}
\begin{proof}
If the assumed conditions hold, then 
\[
\dim_{\C{}} \Cohom{0}{\Xi} = n+1
\]
and the immersion 
\[
E_m \to \CP{}\left(\Cohom{0}{\Xi}^*\right) = \CP{n}
\]
is a local biholomorphism, and under
this map $\mathcal{O}(1)$ pulls
back to $\Xi$. Because $E_m$ is 
compact, it is a covering map. Because
$\CP{n}$ is simply connected, this
map is a global biholomorphism.
\end{proof}

\begin{lemma} \label{lm:Two}
Suppose that $E_m$ is compact. Then
either (1)
every section $\sigma \in \Cohom{0}{\Xi}$
is of the form $\sigma=\sigma_v$
for some $v \in T_m M$ or (2) 
\[
c_1 \left ( \Xi \right )^n > 1
\]
or (3) 
the map
\[
\alpha \in \Cohom{2}{E_m} \to 
\alpha \cap c_1 \left ( \Xi \right )^{n-1} \in
\Cohom{2n-2}{E_m}
\]
has nonempty kernel.
\end{lemma}
\begin{proof}
Suppose that
conditions (1) and (2) do not hold.
Consider a section $\sigma$. Its
zero locus is a projective variety,
since $E_m$ is. Take $P_0 \in E_m$ a smooth
point of
$\left ( \sigma = 0 \right )$.
Pick a $J_{P_0}$ complex basis
$v_0, \dots, v_n$ of $T_m M$ 
with $v_1, \dots, v_n \in P$.
Using the holomorphic immersion 
\[
P \in E_m \mapsto
\left (
\frac{\sigma_{v_1}(P)}{\sigma_{v_0}(P)}, \dots, 
\frac{\sigma_{v_n}(P)}{\sigma_{v_0}(P)}
\right ) 
=
\left ( Z_1, \dots, Z_n \right )
\in \C{n}
\]
(which is only defined in a neighborhood
of $P_0$)
we find that $\left ( \sigma = 0 \right )$
is mapped to an analytic variety
in $\C{n}$, with a smooth point at the
origin. We can make a $J_{P_0}$ complex
linear change of basis to arrange that
$\left ( \sigma = 0 \right )$
is tangent to $\left ( \sigma_{v_1} = 0 \right )$.
Then in the $Z_{\mu}$ coordinates,
\[
\left ( \sigma =0 \right )
=
\left ( Z_1 = f \left ( Z_2, \dots, Z_n \right ) \right ).
\]
The variety
\[
\left ( Z_1 = \dots = Z_{n-1} = 0 \right )
\]
is a straight line, lying entirely 
inside $\left ( Z_1 = 0 \right ).$
So it is tangent to $\left ( Z_1 = f \right )$,
and therefore either lies entirely
inside $\left ( Z_1 = f \right )$
or else strikes it at the origin
with multiplicity (i.e. after small topological
perturbation, strikes at least twice).

By changing the choice of complex basis
$v_0, \dots, v_n$ we can see that the
same is true for any line lying inside
$\left ( Z_1 = 0 \right )$ and
passing through the origin of these
coordinates.

Returning from the $Z_{\mu}$
coordinates, either (1) $\left ( \sigma =  0 \right )$
contains the component of 
$\left ( \sigma_{v_1} = 0 \right )$
passing through $P_0$, or (2) else it has
multiplicity at least two with
$\left ( \sigma_{v_1} = \dots = \sigma_{v_{n-1}} = 0 \right )$
for some choice of $J_{P_0}$ complex basis $v_0,\dots,v_n$.
Topologically, this says that
\[
c_1 \left ( \Xi \right )^n > 1.
\]
But this contradicts our hypotheses.

By positivity of intersections of the
$\left ( \sigma_{v_{\mu}} = 0 \right )$, we see
then that $c_1 \left ( \Xi \right )^n = 1$.
Each $\left ( \sigma_{v_{\mu}} = 0 \right )$
is a smooth variety, since each vector
$v_{\mu}$ belongs to a $J_P$ complex
basis at each point $P \in E_m$.

If $\left ( \sigma_{v_1} = 0 \right )$ has more than
one component, 
say
\[
\left ( \sigma_{v_1} = 0 \right ) = X \cup Y
\]
with $P \in Y$, 
then taking intersection
with the other $\sigma_{v_{\mu}}$ we
find that 
\[
X \cap c_1 \left ( \Xi \right )^{n-1} = 0
\]
contradicting another of our 
topological hypotheses. The same
for $\left(\sigma=0\right)$. So
$\left(\sigma=0\right)=\left(\sigma_{v_1}=0\right)$.
\end{proof}

\begin{corollary}
\label{cor:Teeth}
If $E_m$ is compact,
then either (1) $E_m$ is biholomorphic
to $\CP{n}$ via a biholomorphism
taking $\Xi$ to $\mathcal{O}(1)$,
or (2) 
\[
c_1 \left ( \Xi \right )^n > 1
\]
or (3) the map
\[
\alpha \in \Cohom{2}{E_m} \to 
\alpha \cap c_1 \left (\Xi \right )^{n-1} \in
\Cohom{2n-2}{E_m}
\]
has nonempty kernel.
\end{corollary}

So now we have only to control the
topology of $E_m$ and of $c_1 \left ( \Xi \right )$
and we will be able to control the
complex geometry.

\subsection{A great circle fibration}

Let us continue our study of hypersurface equations
by constructing a great circle fibration on $\Xi^*$.
The line bundle $\Xi \to E_m$ is soldered by
taking any complex hyperplane $P_0 \subset \C{n+1}$,
and forming the quotient of $B_{E_m} \times \C{n+1}/P_0$
by the action of the structure group $H$ on
each of the factors $B_{E_m}$ and $\C{n+1}/P_0$.
By definition of $\Xi$, its fibers are 
\[
\Xi(P) = T_m M/P.
\]
The fibers of $\Xi^*$ are therefore 
\[
\Xi^*(P) = P^{\perp}
\]
where $P^{\perp}$ is the space of 
complex linear 1-forms on $T_m M$ vanishing
on $P$. We see that if we take $P_0 \subset \C{n+1}$
any fixed complex hyperplane, then $\Xi^*$
is the quotient of $B_{E_m} \times P_0^{\perp}$
by the structure group. We have a map
\[
\Phi : B_{E_m} \times P_0^{\perp} \to T'_m M
\]
(where \( T'_m M = \operatorname{Lin}_{\R{}} \left ( T_m M, \R{} \right ) \)
is the real dual space)
defined by
\[
\Phi(\eta,F) = \Real \sum_i F_i \circ \eta.
\]
(where $\Real$ indicates the real part).
We calculate
\[
d \Phi \circ \eta^{-1} = 
\Real
\left ( 
d 
\begin{pmatrix}
F_i
\end{pmatrix}
\right )
- 
\Real
\left ( 
\begin{pmatrix}
F_i &
F_{\bari} 
\end{pmatrix}
\begin{pmatrix}
\omega^i_j & \omega^i_{\barj} & \omega^i_{\mu} & \omega^i_{\barmu} \\
\omega^{\bari}_j & \omega^{\bari}_{\barj} & \omega^{\bari}_{\mu} & \omega^{\bari}_{\barmu} 
\end{pmatrix}
\right )
\]
which, from the structure equations, has full rank
except at $\eta = 0$.
This map descends by $H$ invariance to a map of
$\Xi^* \to T^*_m M$, which is therefore still of full rank.
But then by dimension count, it is a local diffeomorphism.

By writing $S \Xi^*$ I mean the
circle bundle one obtains by looking at
nonzero elements of $\Xi^*$ up to positive real
rescaling. On the other hand, $ST'_m M$
means the sphere constructed
out of $T'_m M \backslash 0$ by quotienting by 
positive real
rescaling. 
\begin{lemma}
If $E_m$ is compact and connected, then
\[
\Phi :  S(\Xi^*) \to  S T'_m M
\]
is a diffeomorphism. 
\end{lemma}
\begin{proof}
The argument above shows that this
map is defined and is a local diffeomorphism.
But if $E_m$ is compact, then so is $S\Xi^*$,
and so the result is immediate since the
sphere $ST'_m M$ is compact and simply connected.
\end{proof}

\begin{lemma} Under the map $\Phi$ the
fibers of $S \Xi^*$ become great circles
on $ST'_m M$. 
\end{lemma}
\begin{proof} Given $\xi \in \Xi^*(P)$
with $\xi \ne 0$, the elements of
the fiber $\Xi^*(P)$ are all of the form
$(a+\sqrt{-1}b)\xi$. They are mapped
to
\[
\Real \left ( a + \sqrt{-1}b \xi \right )
= a \Real \xi - b \Imag \xi
\]
which describes a 2-plane in $T'_m M$.
\end{proof}

Consequently the bundle $S \Xi^* \to E_m$
is a great circle fibration.

\begin{theorem}[McKay \cite{McKay:MoreUnpub}]\label{thm:Mine}
Every smooth great circle fibration of a sphere
is carried by some diffeomorphism to the Hopf
fibration. In particular, the base of the
fibration is diffeomorphic to 
a complex projective space.
\end{theorem}
(With the exception of 3-spheres and 
5-spheres, this result was proven earlier by
C.~T. Yang \cite{Yang:1993}. The result
for 3-spheres is irrelevant for
applications in this article.)

\begin{theorem}
Suppose that $E \subset \Gro{2n}{TM}$
is a system of partial differential equations
for one complex function of several complex
variables, and that $E$ has Cauchy--Riemann tableau.
If the fibers of $E \to M$ are compact and
connected, then $E$ is the Cauchy--Riemann
equations of a unique complex structure
on $M$.
\end{theorem}
\begin{proof}
We have shown now that the bundle $S \Xi^*$
is diffeomorphic to the Hopf fibration. The
bundle $S \Xi^*$ is the principal circle
bundle associated to the complex line bundle
$\Xi^*$. So we must have $\Xi^*$ isomorphic
to the line bundle associated to the Hopf
fibration, which is the line bundle $\mathcal{O}(-1)$.
Therefore $\Xi$ is isomorphic to $\mathcal{O}(1)$,
as a complex vector bundle on a real manifold,
and its Chern number is
\[
c_1 \left ( \Xi \right ) = 1.
\]
By corollary~\vref{cor:Teeth},
we find that $E_m$ is biholomorphic
to $\CP{n}$ via a biholomorphism
taking $\Xi$ to $\mathcal{O}(1)$.
The rest follows from theorem~\vref{thm:BIG}.
\end{proof}

\begin{corollary} Every Legendre fibration
$E \to M$
(not assumed to be holomorphic)
with compact connected fibers 
on a complex contact manifold $E$
provides $M$ with
a unique complex structure for which
there is a holomorphic contactomorphism
\[
\xymatrix{
E \ar[rr] \ar[dr] & & J^1 M \ar[dl] \\
                  & M &  
}
\] 
where $J^1 M$ is the bundle of 
projectivized holomorphic cotangent
spaces of $M$. In particular, the
fibers of $E \to M$ are complex projective
spaces.
\end{corollary}

\subsection{A Chern class}

We will now present yet another approach to proving
that there are no hypersurface
equations with compact fibers.
This gets around the use of
Merkulov's results on complex
contact geometry, but still
requires our theorem
on great circle fibrations.

\begin{definition}
Suppose that $E \subset \Gro{2n}{TM}$ is
a system of partial differential equations
with Cauchy--Riemann tableau, for one
complex function of several complex
variables. We will say that $E$ has 
\emph{Cauchy--Riemann local
topology} if each fiber $E_m$ is compact and 
on each of the
fibers $E_m \subset \Gro{2n}{T_m M}$
the complex line bundle $\Xi=TE/\Theta$ satisfies
\[
(n+1) c_1 \left ( \Xi \right ) + c_1 \left ( K \right ) \le 0
\]
where $K$ is the canonical bundle of $E_m$.
\end{definition}

This condition is truly topological, and does
not require any information about the biholomorphism
type of the fiber. Also, if it holds at one
fiber, and all of the fibers are compact,
and the base $M$ is connected, then
it holds at all fibers.

\begin{proposition} \label{prop:BigTwo}
Suppose that $E \subset \Gro{2n}{TM}$ is
a system of partial differential equations
with Cauchy--Riemann tableau,
for one complex function of several complex
variables. Suppose that $E$ has
Cauchy--Riemann local topology.
Then there is a unique complex
structure on $M$ for which $E$ is the
Cauchy--Riemann equations for complex
hypersurfaces. 
\end{proposition}
\begin{proof}
Consider the line bundle over 
$E_m$ soldered by the expression $\det a^i_j \det a^{\mu}_{\nu}$.
This is the line bundle whose
fiber over a point $P \in E_m$ consists of
$\Det_J V$ where the $\Det_J V$ is the complex
determinant for the complex structure $J=J(P)$ on
$V$ determined by that point of $E_m$.
Then the 1-form
\[
A = \omega^i_i + \omega^{\mu}_{\mu}
\]
is a connection 1-form for that line
bundle, and its curvature 2-form
is
\begin{align*}
F =& \frac{i}{2 \pi} dA \\
=& \frac{i}{2 \pi} t^{\mu \epsilon}_{\barsigma m} t^{\barsigma \bartau}_{\mu \barp} 
\omega^m_{\epsilon} \wedge \omega^{\barp}_{\bartau} \\
\ge & 0.
\end{align*}
Therefore the line bundle $\Det_J$ has
nonnegative Chern classes.
Moreover if it has vanishing first Chern class, then
we must have $F=0$, so must have $t^{\mu \epsilon}_{\barsigma m}=0$.
From the structure equations we see that all of the
invariants of $E$ vanish, 
so by the Newlander--Nirenberg theorem and Darboux's theorem,
% or Malgrange's theorems in 
% Malgrange
% \cite{Malgrange:1972a,Malgrange:1972b},
we can easily see that $E$ is the Cauchy--Riemann
equations.

Now we have only to ascertain the relation
between the Chern classes of $\Det_J V$
and those of $\Xi$. We leave
it to the reader to show that
\[
\Det_J V = \Xi^{\otimes (n+1)} \otimes_{\C{}} K
\]
as complex line bundles (evident from the soldering
representations).
\end{proof}

\begin{corollary} Any hypersurface equation 
$E \to \Gro{2n}{TM}$ with compact fibers 
must be the Cauchy--Riemann equations
of a unique complex structure on $M$.
\end{corollary}
\begin{proof} First we apply
the great circle fibration theorem~\vref{thm:Mine} 
to identify 
the fiber $E_m$ with
$\CP{n}$ diffeomorphically,
and identify $\Xi$ with
$\mathcal{O}(1)$ diffeomorphically.
This determines the topology
of $\Xi$ and $E_m$ completely,
and allows us to apply 
corollary~\vref{cor:Teeth}
to see that via biholomorphism,
we can identify $E_m$ with $\CP{n}$ and
$\Xi$ with $\mathcal{O}(1)$. 
This determines the Chern 
class of the canonical 
bundle, allowing us to
employ proposition~\vref{prop:BigTwo}.
\end{proof}

In particular, the moduli space of
hypersurface equations on $M$ with compact
fibers is canonically globally isomorphic to the
moduli space of complex structures on 
$M$.

\section{Moduli of hypersurface equations}
There are no particularly interesting
hypersurface equations, due to the absence
of compactness in the fibers. Nonetheless,
we may ask how many such equations exist,
and whether we can deform them into one
another.

A first approach to this question: give a
complex contact manifold $E$, we can
ask how many fiber bundle mappings $E \to M$ 
it admits which determine hypersurface
equations (on the base $M$). This is just
asking for not-necessarily-holomorphic
Legendre fibrations by holomorphic Legendre
submanifolds. There is little one can say
about the global problem, other than that
the fibers can not be compact without
entering into the territory of our theorems
above. But the local study is straightforward.
In local coordinates $w,z_{\mu},p_{\mu}$
in which the contact structure is
$dw - p_{\mu} \, dz^{\mu}=0$, every
Legendre submanifold transverse to 
$w=z=0$ is of the form $w=w(z), p_{\mu}=\pd{w}{z^{\mu}}$.
We can arrange such transversality by
change of coordinates, and even by
arbitrarily small linear change of coordinates.
A Legendre fibration will have a unique
Legendre submanifold passing through 
each point $(w,z,p)=(0,0,P)$. So we will
have a function $w\left(z,P,\bar{P}\right)$
for which $w\left(0,P,\bar{P}\right)=0$
and $\pd{w}{z^{\mu}}\left(0,P,\bar{P}\right)=P_{\mu}$.
Hence 
\[
\pd{}{z^{\nu}}\left( w - P_{\mu} z^{\mu} \right) = 0,
\]
so 
\[
w = P_{\mu} z^{\mu} + f_{\mu \nu} \left(z,P,\bar{P}\right)
z^{\mu} z^{\nu},
\]
for some functions $f_{\mu \nu}\left(z,P,\bar{P}\right)$,
holomorphic in $z$ and smooth in $P$.
These functions $f_{\mu \nu}$ can be chosen arbitrarily.
Clearly the
space of germs of Legendre fibrations is connected.
So there is a connected moduli space of germs
of hypersurface equations, from this point of view.
The local invariants described above show that the
moduli space is not a single point.
Since the automorphisms of $E$ must be complex analytic,
while the $f_{\mu \nu}$ need not be, the moduli space
is of infinite dimension.

Another point of view, along the lines of
writing down partial differential equations,
would ask for hypersurface equations on a fixed
manifold $M$. 
A hypersurface equation is a submanifold
$E \subset \Gro{2n}{TM}$, satisfying
an overdetermined first-order system of equations.
Again, if we are willing to allow
noncompact fibers, then a global characterisation
seems impossible, and we restrict attention
to germs. Clearly the constructions above
impose an orientation on the base manifold $M$.
So far, this is the only obstruction we
have found to prevent deformation of
one hypersurface equation into another.
We conjecture that this is the only 
deformation obstruction for hypersurface
equation germs.

\section{Directions for further research}
The next obvious step is to investigate
the analogues of calibrated submanifolds
of exceptional holonomy manifolds 
in the same sense. Calibrated submanifolds
are submanifolds of a Riemannian manifold
which are not only minimal submanifolds, but
satisfy a first-order system of partial
differential equations which forces them
to be of minimal volume globally, among
all submanifolds in the same homology class
(or perhaps relative homology class, fixing
boundary components).
The relevant first-order equations are more flexible
than the Cauchy--Riemann equations.
For example, in Calabi--Yau manifolds, the
special Lagrangian submanifolds are calibrated.
One can already deform (at least locally)
the Calabi--Yau metric, and so obtain 
infinite dimensional families of equations
for calibrated submanifolds, all with the same
tableau. But the true flexibility of calibrated
submanifold equations is unknown. If the
equations are flexible, with fixed tableau,
this may hold the key to large perturbation
theorems. The point is to 
deform the special Lagrangian equation, 
preserving tableau, which could
be much easier than deforming the
Calabi--Yau metric.

Finally, the theorems proven here give 
impetus to the development of the analogy between
the theory of complex surfaces and the theory of
elliptic systems of partial differential equations
for 2 functions of 2 variables.
Complex surfaces are deeply understood.
By contrast, 
these elliptic equations have unusually easy analysis,
but few global theorems. The use of these ideas
(1) to classify the diffeomorphism types of smooth projective
planes (see McKay \cite{McKay:2005}) and (2) in 
Gromov's nonlinear Riemann
mapping theorem \cite{Gromov:1985}, p. 144, 
shows that this
analogy is potent but still in its infancy. 

% \nocite{*}
\bibliographystyle{amsplain}
\bibliography{higher}

\end{document}